\documentclass[preprint,12pt]{elsarticle}
\usepackage{amssymb,amsmath}
\usepackage{fullpage,graphicx,color,float}
\usepackage[dvipsnames]{xcolor}
\usepackage[percent]{overpic}
\numberwithin{equation}{section}

\newcommand\alf{{\alpha}}
\newcommand\alfh{\widehat{\alpha}}
\newcommand\betah{\widehat{\beta}}
\newcommand\alfb{\alpha_b}
\newcommand\betb{\beta_b}
\newcommand\calf{c_{\alpha}}
\newcommand\cbet{c_{\beta}}
\newcommand\bet{{\beta}}
\newcommand\alfp{\alf_b}
\newcommand\betp{\bet_b}
\newcommand\xb{\mathbf{x}}
\newcommand\nw{n_w}
\newcommand\xbtil{\tilde{\mathbf{x}}}
\newcommand\fbtil{\tilde{\mathbf{f}}}
\newcommand\xotil{\tilde{\mathbf{x}}_0}
\newcommand\rhotil{\tilde{\rho}}
\newcommand\ktil{\tilde{k}}

\newcommand\Xb{\mathbf{X}}
\newcommand\Ab{\mathbf{A}}
\newcommand\Rb{\mathbf{R}}

\newcommand\sbk{\mathbf{s}_k}
\newcommand\xo{\xb_0}
\newcommand\xbh{\widehat{\xb}}
\newcommand\xbhtil{\widehat{\xbtil}}

\newcommand\rhoo{\rho_0}
\newcommand\xp{\xb_b}
\newcommand\xbb{\xb_b}
\newcommand\xalf{\xb_{\alf}}
\newcommand\xbet{\xb_{\bet}}

\newcommand\nb{\mathbf{n}}
\newcommand\fb{\mathbf{f}}
\newcommand\ub{\mathbf{u}}
\newcommand\yb{\mathbf{y}}

\newcommand\Ubinf{\mathbf{u}_{\infty}}

\newcommand\delalf{\Delta\alpha}
\newcommand\delbet{\Delta\beta}
\newcommand\nwin{n_w}
\newcommand\sumtrap[2]{\sum_{#1=0}^{#2}\hskip-5pt\hbox{ }'}

\newcommand\kmax{K}
\newcommand\mmax{M}

\newcommand\etal{et al.\ }

\begin{document}
\begin{frontmatter}

\title{Corrected Trapezoidal Rules for Near-Singular Surface Integrals Applied to 3D Ellipsoids in Stokes Flow}

\author{Monika Nitsche\corref{cor1}\fnref{label1}}
\cortext[cor1]{corresponding author email: nitsche@unm.edu}
\affiliation[label1]{organization={Department of Mathematics and Statistics},
           addressline={University of New Mexico}, 
           city={Albuquerque},
           postcode={87108}, 
           state={NM},
           country={USA}}
           
\author{Bowei Wu\fnref{label2}}
\affiliation[label2]{organization={Department of Mathematics and Statistics},
           addressline={University of Massachusetts Lowell}, 
           city={Lowell},
           postcode={01854}, 
           state={MA},
           country={USA}}

\author{Ling Xu\fnref{label3}}
\affiliation[label3]{organization={Department of Mathematics and Statistics},
           addressline={North Carolina A \& T State University}, 
           city={Greensboro},
           postcode={27411}, 
           state={NC},
           country={USA}}

\begin{abstract}
Interfacial Stokes flow can be efficiently computed using the Boundary Integral Equation method. In 3D, the fluid velocity at a target point is given by a 2D surface integral over all interfaces, thus reducing the dimension of the problem. 
A core challenge is that for target points near, but not on, an interface, the surface integral is near-singular and standard quadratures lose accuracy.  This paper presents a method to accurately compute the near-singular integrals arising in elliptic boundary value problems in 3D.  It is based on a local series approximation of the integrand about a base point on the surface, obtained by  orthogonal projection of the target point onto the surface. The elementary functions in the resulting series approximation can be integrated to high accuracy in a neighborhood of the base point using a recursive algorithm. The remaining integral is evaluated numerically using a standard quadrature rule, chosen here to be the 4th order Trapezoidal rule. The method is reduced to the standard quadrature plus a correction, and is uniformly of 4th order. The method is applied to resolve Stokes flow past several ellipsoidal rigid bodies. We compare the error in the velocity near the bodies, and in the time  and displacement of particles traveling around the bodies, computed with and without the corrections.
\end{abstract}

\end{frontmatter}

\section{Introduction}
The Stokes equation modeling highly viscous flow in the zero-inertia limit is an elliptic partial differential equation whose solution can be computed using the Boundary Integral Equation (BIE) method \cite{kress2014linear,pozrikidis1992boundary}.  The BIE approach is especially efficient for particulate fluid flows -- viscous flows that contain soft or rigid bodies such as red blood cells, vesicles, drops, bubbles, colloidal particles, and micro swimmers -- in which constraints placed on the interface between body interior and exterior enter the boundary integral representation of the fluid velocity. Accurate resolution of near fluid-particle or particle-particle interactions is essential for the accuracy of the overall numerical scheme and the stability of long-term simulations. A core challenge is the accurate evaluation of the Stokes boundary integral operators, also known as layer potentials, at target points close to, but not on, the boundary surface. These integrals are near-singular and notoriously difficult to approximate  using standard quadrature rules.

Considerable research efforts have focused on developing efficient near-singular quadrature schemes for integral operators associated with elliptic PDEs. 
These quadrature methods can be loosely categorized into four classes.
The first class of methods exploits the smoothness of layer potentials in the far-field. Kl\"ockner \etal \cite{klockner2013} introduced the QBX method {in 2D, later generalized to 3D \cite{afKlinteberg2016,wala2020}}, in which an accurate series expansion is constructed in the far field and then evaluated at target points near or on the boundary. Similarly, the piecewise surface quadratures of Ying \etal \cite{ying2006high} and of Morse \etal \cite{morse2021robust} use interpolation or extrapolation of the far-field integral values to the near-field.
A second class of methods is based on harmonic  analysis. 
{In two dimensions, Helsing and Ojala \cite{helsing-ojala2008} developed specialized Gauss-Legendre panel quadratures} using complex interpolation to write near-singular line integrals as a series of integrals that can be evaluated analytically using recursion. {Later, the Singularity Swap Quadrature (SSQ) developed by af Klinteberg and Barnett \cite{afKlinteberg2020} improved the convergence of the Helsing-Ojala quadrature on curved panels by a singularity cancellation in the complexified parameter space; the SSQ also applies to line integrals in 3D and has been extended to the trapezoidal rule \cite{afKlinteberg2024} and to axisymmetric surfaces \cite{krantz_tornberg}.  On the other hand, Zhu and Veerapaneni \cite{zhu2022quaternion} proposed an innovative extension of the Helsing-Ojala quadrature} to 3D triangulated surfaces using harmonic polynomial interpolation based on quaternionic algebra and analysis. Perez-Arancibia \etal \cite{perez-arancibia2019interpolation} used  interpolation of the density function together with Green's identities to regularize the integrand. The global quadrature of Barnett \etal \cite{lsc2d} expresses the layer potentials as Cauchy integrals and evaluates them using a spectrally accurate barycentric formula derived from singularity subtraction.
A third class of methods are the regularized quadratures of Beale and collaborators \cite{beale2001method,beale2024extrapolated,beale_tlupova} who compute regularized singular or near-singular integrands and then correct the regularization errors using asymptotic approximations. The method presented in this paper belongs to a fourth class: error corrections for regular quadratures applied to singular or near-singular integrals. Carvalho \etal \cite{carvalho2018,khatri2020close} find asymptotic error expansions for near-singular integrals on closed curves in 2D based on the periodic trapezoidal rule, and on 3D surfaces based on a local auxiliary polar mesh. Most other existing error correction methods are based on the trapezoidal  rule and are for singular integrals \cite{duan2009high,wu2021zeta,wu2021corrected,izzo2022trap}. Wu \cite{wu2024extension} has recently generalized the Euler-Maclaurin formulas for singular integrals \cite{sidi2012algebraic} to near-singular integrals, using a connection between polygamma functions and the Riemann zeta series. Closely related are the methods of Nitsche \cite{nitsche2021tcfd,nitsche2022acom}, where the near-singular errors are obtained using local approximations of the integrand that can be integrated exactly. With this approach, corrected trapezoidal rules were obtained for planar vortex sheet flow and for axi-symmetric Stokes flow. They were recently applied  to compute the motion of interacting flapping plates \cite{nitscheozasiegel2024stability} and of freely falling plates \cite{sohn2024fallingtcfd}. 

\begin{figure}
 \centering
\vskip -0.2truein
\includegraphics[trim=0 0 30 0, clip, width=0.38\textwidth]{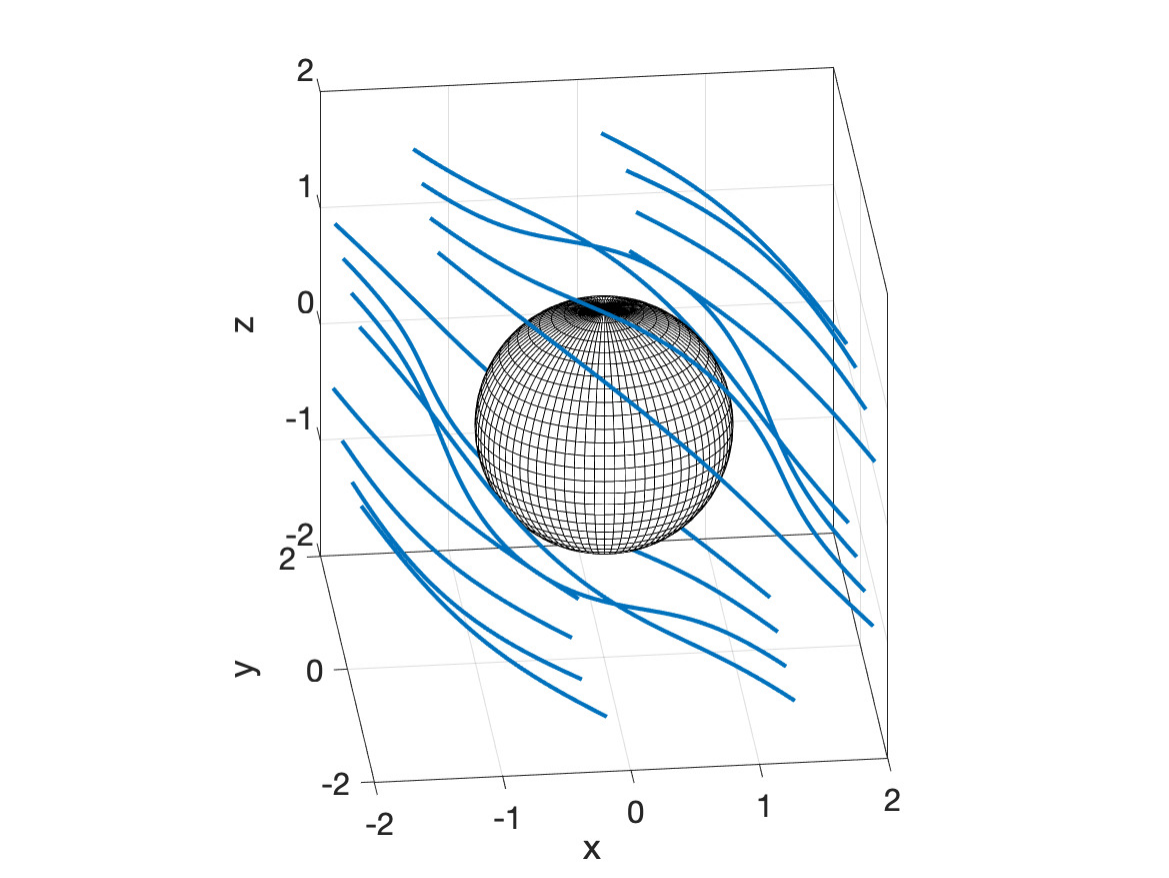}
\includegraphics[trim=28 02 50 50, clip, width=0.39\textwidth]{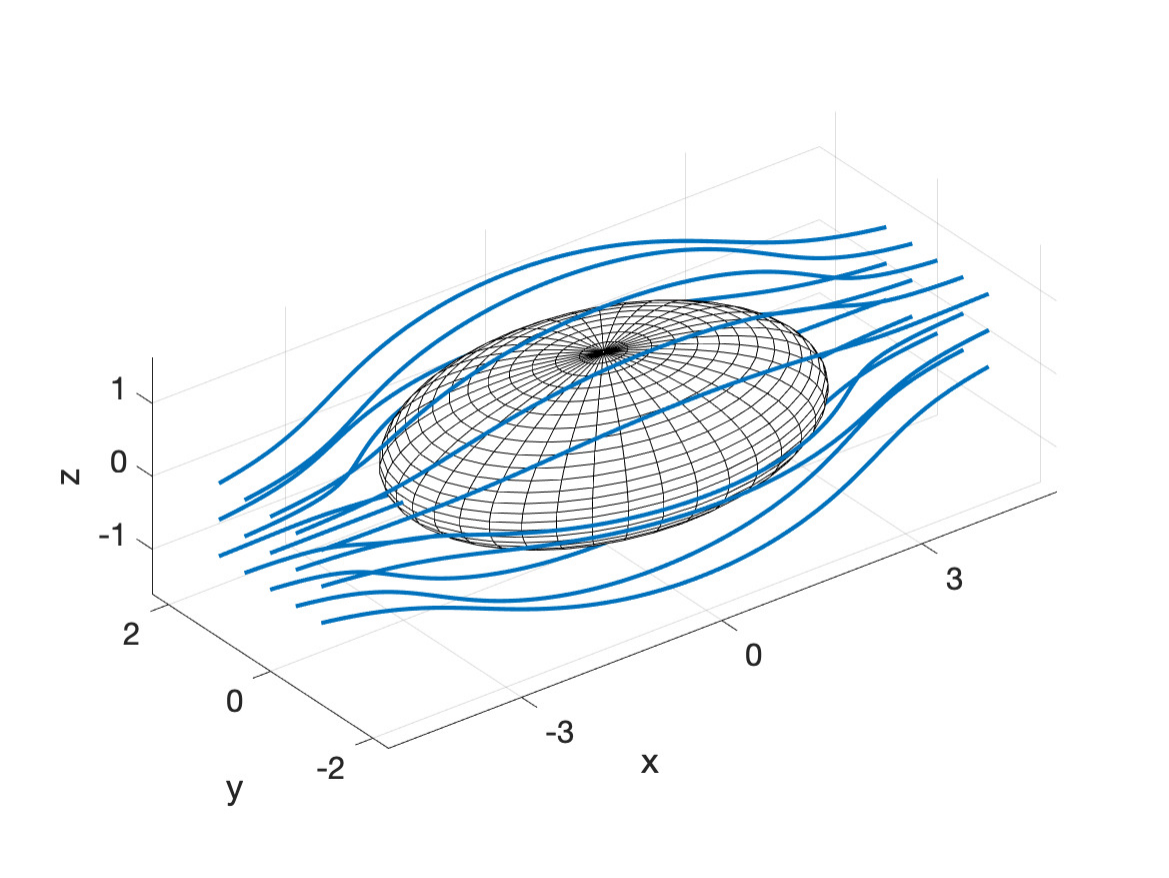}
\includegraphics[trim=68 32 80 30, clip, width=0.34\textwidth]{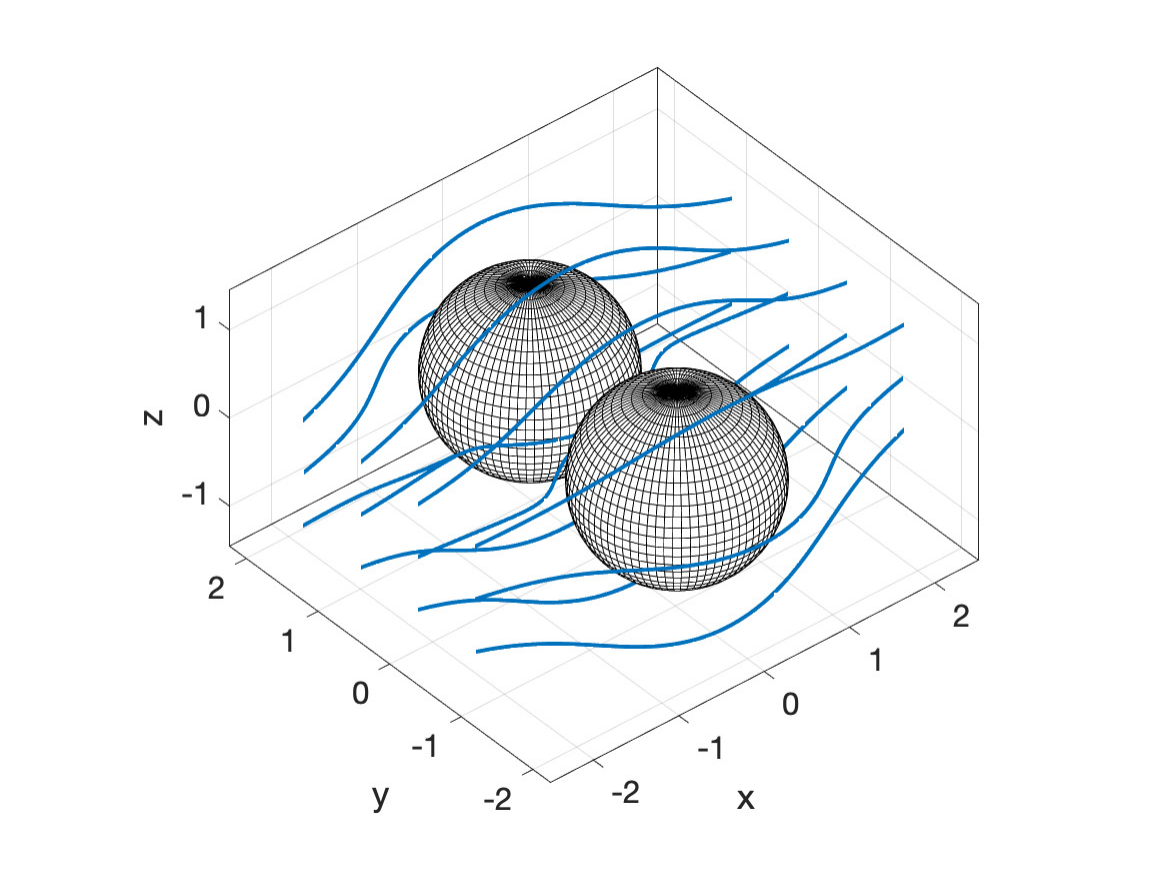}
\includegraphics[trim=33 20 45 46, clip, width=0.45\textwidth]{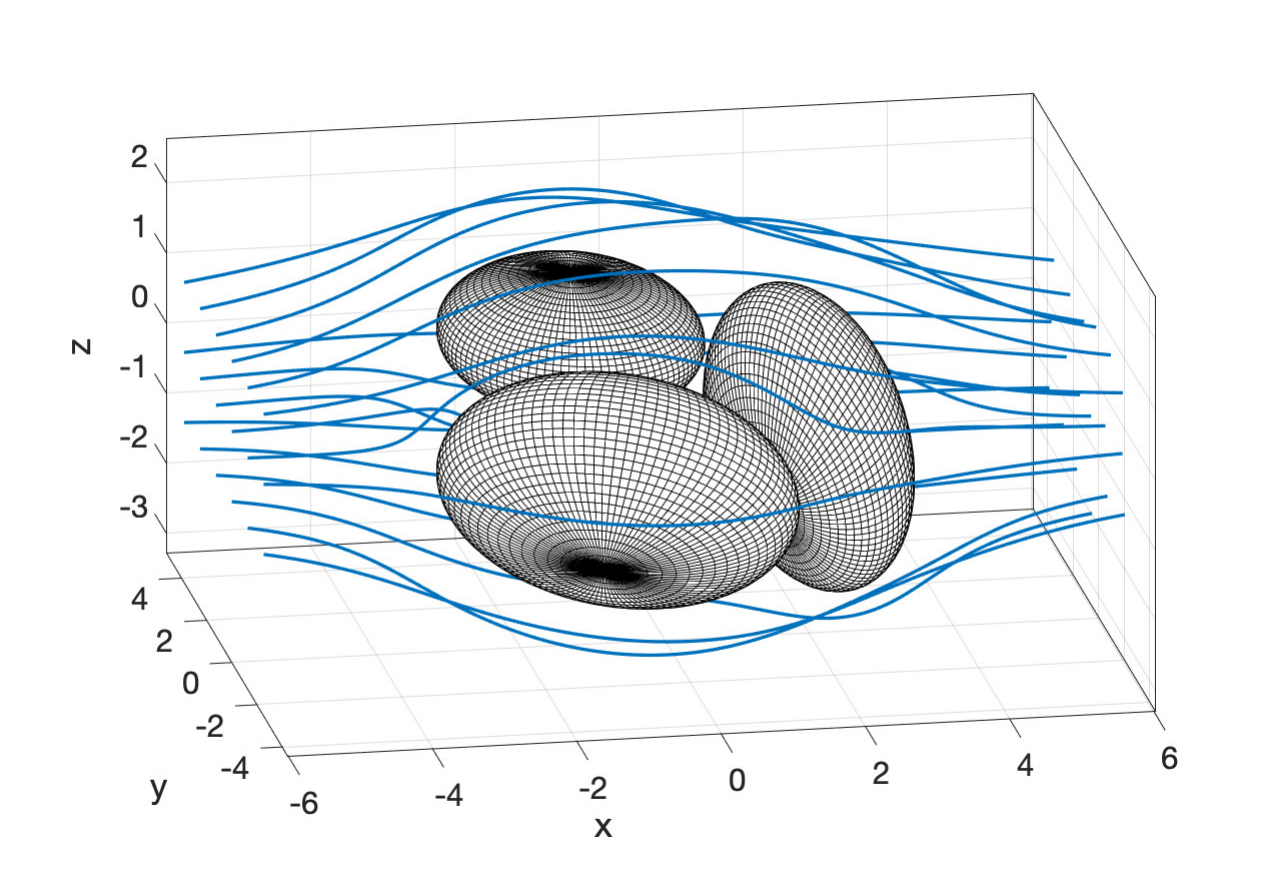}
\vskip-3.8truein \hbox{ }\hskip-1.7truein(a)\hskip2.2truein(b)
\vskip1.4truein \hbox{ }\hskip-2.6truein(c)\hskip1.6truein(d)
\vskip1.6truein
\label{F:examples}
\caption{
Stokes flow (a) past sphere, (b) past ellipsoid, (c) past  two spheres, (d) past three ellipsoids.}
\end{figure}

This paper generalizes the latter approach to evaluate surface integrals in 3D of the form
\begin{equation}
\int_S\frac{F(\xb,\xo)}{|\xb-\xo|^r}\,\omega(\xb)dS\,,
\label{E:genform}
\end{equation}
for odd integers $r=1,3,5,\dots$ that occur in the BIE approach for elliptic problems. These integrals are near-singular when the target point $\xo$ is near, but not on, $S$. They are integrable in the principal value sense when $\xo$ is on $S$. The method consists of using Taylor series to find local approximations of the integrand about the orthogonal projection of the target point onto the surface. The basis functions in the resulting expansion capture the near-singularity of the integrand and are integrated {efficiently and} to high accuracy in a neighbourhood of the projection. The remaining integral is computed using a standard quadrature, chosen here to be a 4th order trapezoidal rule. The number of terms kept in the Taylor approximation of the integrand is determined using upper bounds on the trapezoidal errors of their integrals. The result reduces to the standard trapezoidal rule plus a correction, and is uniformly of 4th order. The method is applied to resolve Stokes flow past several ellipsoidal rigid bodies, for the four examples illustrated in figure \ref{F:examples}. We confirm 4th order convergence in the fluid velocity near the bodies, and determine the effect of the corrections on the time  and displacement of particles traveling around the bodies. {Note that although we only present rigid ellipsoids in these numerical examples, our quadrature method can also apply without modifications to broader examples such as deformable particles in Stokes flow \cite{veerapaneni2011}, which will be considered in future work.}

As mentioned, the method introduced here generalizes that for line integrals introduced in \cite{nitsche2021tcfd}. In that case, the basis  functions 
of the local approximations can be  integrated exactly. In the present case of surface integrals, the basis functions do not have readily available antiderivatives. Inspired by the recurrence relations for singular integrals in Helsing \cite{helsing2013higher}, we derive a recurrence relation for a family of near-singular integrals that returns highly accurate integral values. Details of the derivation will be presented in forthcoming work. 

The paper is organized as follows. Section 2 formulates the fluid velocity in terms of single and double layer potentials. Section 3 presents the numerical method, including  the local approximation, the recurrence relation used to integrate the basis function, and details of the implementation. Section 4 presents numerical results for the four sample Stokes flows. A summary is presented in Section 5.

\section{Stokes Flow Past Objects}
Consider Stokes flow past a body with surface $S$, driven by a far-field flow $\ub_{\infty}$. The velocity at a \textit{target point} $\xo$ is given as a perturbation from $\ub_{\infty}$ induced by the body, where the perturbed flow is given by a compound representation \cite[Eqs 4.7.1 and 4.7.4]{pozrikidis1992boundary} of single and double layer potentials (SLP and DLP), 
\begin{equation}
\ub(\xo)= \ub_{\infty}+{\cal{S}}[\fb](\xo) +{\cal{D}}[\fb](\xo) \,,
\label{E:veloone}
\end{equation}
where
\begin{subequations}
\begin{align}
{\cal{S}}[\fb](\xo)&=\frac{1}{8\pi\mu}\int_S \Big(\frac{\fb}{\rho}+\frac{(\fb\cdot\xbh)\,\xbh }{\rho^3}\Big) \,dS(\xb)~,\quad\\[2ex]
{\cal{D}}[\fb](\xo)&=-\frac{3}{4\pi}\int_S\frac{(\fb\cdot\xbh)\,\xbh\,(\xbh\cdot\nb)}{\rho^5}\,dS(\xb)~.
\label{E:layerpotentials}
\end{align}
\end{subequations}
Here $\rho^2=|\xb-\xo|^2$ and $\xbh=\xb-\xo$. That is, the integrals are of the form (\ref{E:genform}) with $r=1,3,5$, and are near-singular if $\xo$ is near $S$. The compound formulation \eqref{E:veloone} will lead to the well-posed boundary integral equation \eqref{E:boundmult}, where the use of the SLP helps eliminate the 6-dimensional nullspace of the DLP corresponding to rigid body motions of $S$ (see discussions in \cite[p.129]{pozrikidis1992boundary}).

For flow past several objects $E_k$, $k=1,\dots,N$,  such as illustrated in figure \ref{F:examples}(c,d), the velocity is
\begin{equation}
\ub(\xo)= \ub_{\infty}+\sum_{k=1}^N\Big({\cal{S}}_k[\fb_k](\xo) +{\cal{D}}_k[\fb_k](\xo) \Big)
\label{E:velomult}
\end{equation}
and ${\cal{S}}_k$ and ${\cal{D}}_k$ are the Stokes SLP and DLP on the object $E_k$. 
In all our examples below, the far field $\ub_{\infty}$ is constant and the body velocity is zero.  

The functions $\fb_k$ are referred to as the density on body $E_k$.
They are determined by the no-slip boundary condition: the fluid velocity equals the body velocity on all $S_k$. For zero body velocity,
\begin{equation}
\ub(\xb)=\mathbf{0} \hbox{ for all } \xb\in S_k\,,\quad k=1,\dots,N.
\label{E:boundary}
\end{equation}
For flow past a sphere, the density is known analytically to be $\fb=-1.5\ub_{\infty}$ \cite{ockendon1995}.
In all other cases, $\fb_k$ is obtained numerically as follows. 
{Each $S_k$ is discretized on a standard spherical grid $\Xb_k$ of order $p$, which consists of $(2p+2)$ equally spaced points along the azimuthal coordinate and $(p+1)$ Gauss-Legendre points in (the cosine of) the polar angle coordinate. Let $\xb_j, j=1,\dots,P$, denote the $P:=(p+1)(2p+2)$ grid points on $\Xb_k$. The integrals $\cal S$ and $\cal D$ are discretized and collocated on the same spherical grids $\Xb_k$, then
(\ref{E:boundary}) reduces to a linear system of $N\times N$ blocks, each block of size $P\times P$, for the unknowns $\fb_k$,}
\begin{subequations}
\begin{equation}
\frac{1}{2}\fb_k + 
\sum_{k'=1}^N \Ab_{kk'}\fb_{k'}
= -\ub_{\infty}~,\quad k=1,\dots,N\,,
\label{E:boundmult}
\end{equation}
where 
\begin{equation}
\Ab_{kk'}\fb_{k'}={\cal{S}}_{k'}[\fb_{k'}](\xb_j) +{\cal{D}}_{k'}[\fb_{k'}](\xb_j)~,\quad 
\xb_j\in \Xb_k\,.
\end{equation}
\end{subequations}
The term $\frac{1}{2}\fb_k$ in (\ref{E:boundmult}) accounts for the jump condition in the double layer when $\xb_j$ lies on $\Xb_k$ \cite[Eq 2.1.12]{pozrikidis1992boundary}. The diagonal interactions $\Ab_{kk}\fb_k$ describe the self-induced velocity of the body $E_k$ at a point $\xb_j$ on itself, and the off-diagonal interactions $\Ab_{kk'}\fb_{k'}$, $k\ne k'$, describe the velocity induced by $S_{k'}$ at points on $S_k$. The system (\ref{E:boundmult}) is solved by GMRES iterations. In each iteration,  the singular integrals $\Ab_{kk}\fb_k$ are discretized using the spectrally accurate quadrature of \cite{graham2002} and evaluated using the fast grid rotation algorithm of \cite{gimbutasveerapaneni2013rotations} which is based on the Non-uniform Fast Fourier Transform (NUFFT) \cite{greengard2004} and has a time complexity of $O(P^2\log P)$ per body. The off-diagonal interactions $\Ab_{kk'}\fb_{k'}$, $k\ne k'$, are discretized by standard spherical quadrature if the bodies $E_{k}$ and $E_{k'}$ are far apart. If $E_{k}$ is sufficiently close to $E_{k'}$, then $\Ab_{kk'}\fb_{k'}$ is nearly singular for any $\xb_j\in \Xb_k$ that is near $E_{k'}$, which must be evaluated using the corrected trapezoidal quadrature method described in the next section; this trapezoidal quadrature requires first interpolating the density $\fb_{k'}$ onto one of the two uniform grids \eqref{E:grids} on $S_{k'}$, which incurs an additional $O(P\log P)$ cost via NUFFT on each body $E_{k'}$ (independent of the number of points close to $E_{k'}$) -- a negligible cost compared to evaluating the diagonal interactions.

\section{Numerical Method}
\subsection{Basic Overview}
Consider an integral of form (\ref{E:genform}) over a body $E$ with boundary surface $S$. We assume that $S$ is $logically$ $rectangular$, that is, it can be parametrized by $\xb(\alf,\bet)$ where $(\alf,\beta)$ range over a rectangle $D$. After the change of variables, the integral over $D$ is approximated by a  quadrature $T_D$, chosen here to be a high-order trapezoidal rule, 
\begin{equation}
\int_D \frac{F(\xb(\alf,\bet),\xo)}{\rho^r}\omega(\xb(\alf,\bet))J(\alf,\bet)\,d\alf\,d\bet 
=\int_D G(\alf,\bet)\,d\alf\,d\bet\approx T_D[G]~.
\label{E:changeofvars}
\end{equation}
Here $J=|\xb_{\alf}\times\xb_{\bet}|$, $G$ denotes the integrand, and the subscript on $T_D$ denotes the domain of integration.

If the target point $\xo$ is far from the surface $S$, the chosen trapezoidal rule yields uniformly high accuracy. If the target point is so close that the integral is near singular, the trapezoidal rule loses accuracy. In that case, we find a local approximation $H$ of $G$ that (1) captures the near-singularity and (2) can essentially be integrated exactly over a subdomain (or window) $W$ of $D$.  The smoother function $G-H$ is then integrated numerically over $W$, while $G$ is integrated numerically over the remaining domain.  This is summarized as 
\begin{equation}
\begin{split}
\int\displaylimits_D\hskip-3pt G
&=\int\displaylimits_{D\setminus W}\hskip-6pt G
\, +\int\displaylimits_{W}\hskip-3pt (G-H)
\, +\int\displaylimits_{W}\hskip-3pt H\\
&\approx T_{\scriptsize D\setminus W}\hskip-1pt[G]
+T_{\scriptsize W} [G-H]
+\int\displaylimits_{W}\hskip-3pt H\\
&= T_{\scriptsize D}[G]+E_{\scriptsize W}[H],
\label{E:method}
\end{split}
\end{equation}
where
\begin{equation}
E_W[H]=\int\displaylimits_W H-T_W[H].
\end{equation}
That is, it amounts to simply computing the trapezoid approximation $T_D[G]$ and adding a correction $E_{\scriptsize W}[H]$. The window $W$ is  an approximately square subgrid of $D$. For $T_D[G]$ we use a 4th order trapezoidal rule (given  in \S 3.5). For reasons explained in \S 3.6, we use a higher 6th order approximation for $T_W[H]$. The final results are uniformly 4th order accurate over all target points $\xo$. {In practice, the correction is applied if the distance $d$ from the target point to the surface $S$ is small relative to the physical grid size $\Delta s$ of the surface near the target point (see \S 3.8).}

The approximation $H$ is obtained by expanding $G$ about the base point $(\alfb,\betb)$ determined by the orthogonal projection of $\xo$ onto $S$, using Taylor series. The resulting expansion
\begin{subequations}
\begin{equation}
G(\alf,\bet)\approx H(\alf,\bet)=\sum c_{pqk}H_{pqk}(\alf-\alf_b,\bet-\bet_b)~,
\end{equation}
is given
in terms of basis functions 
\begin{equation}
H_{pqk}(\alf,\bet)=\frac{\alf^p\bet^q}{\rho_o^{2k+1}}~,\quad
\rho^2_o=d^2 +c^2_{\alf}\alf^2 +2c_{\alf\bet}\alf\bet +c^2_{\bet}\bet^2,
\label{E:basisfunctions}
\end{equation}
\label{E:expansionG}%
\end{subequations}
where $d$ is the distance from $\xo$ to $S$. These functions capture the near-singular behaviour of $G$. Moreover, while analytic expressions for their  antiderivatives are not readily available, a recursion for $\int_W H_{pqk}$ enables quick and accurate integration. All necessary details are given next. 

\subsection{Parametrization and discretization}
\begin{figure}
 \centering
\includegraphics[trim=0 0 0 0, clip, width=0.85\textwidth]{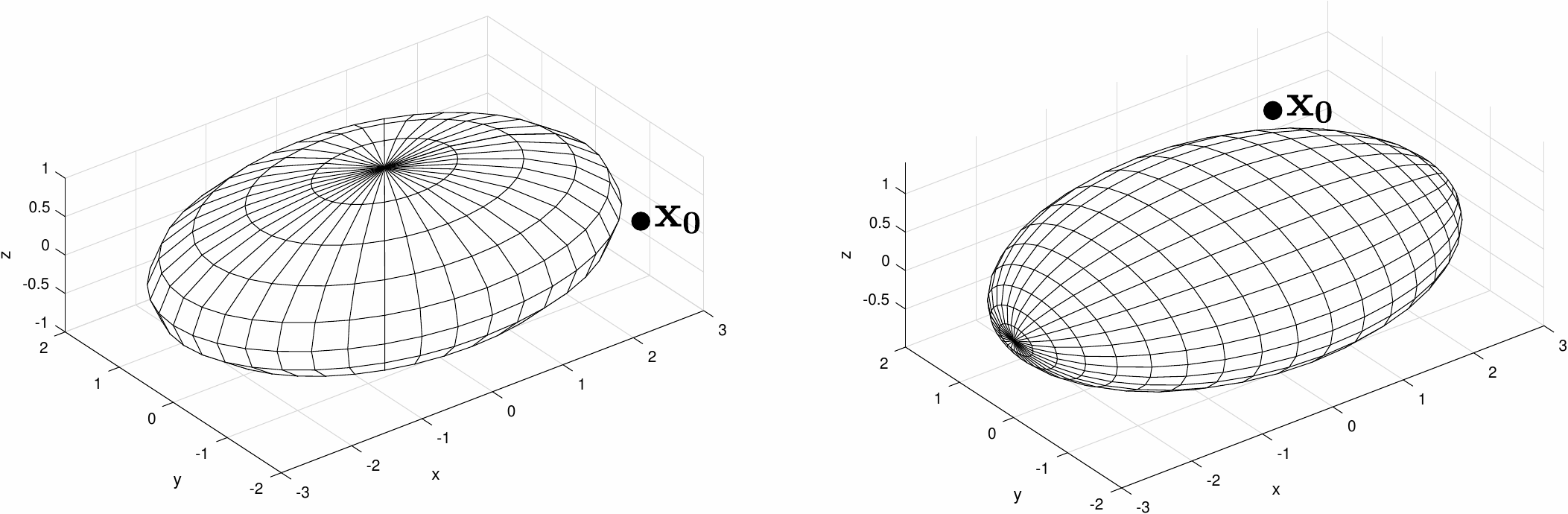}
\vskip-1.8truein \hbox{ }\hskip-1.8truein(a)\hskip2.8truein (b)
\vskip1.4truein
\caption{Sample set of two grids for a standard ellipsoid with $a=3$, $b=2$, $c=1$, and sample target point $\xo$.
(a) Grid 1, with $n_1=10$, $m_1=40$, poles on $z$-axis.
(b) Grid 2, with $n_2=20$, $m_2=30$, poles on $x$-axis.
For given $\xo$, choose the grid whose poles are furthest from $\xo$.}
\label{F:twogrids}
\end{figure}

\paragraph{Standard ellipsoid} For simplicity, all objects considered here are ellipsoids.  Consider first a single ellipsoid $E$ centered at the origin of a Cartesian coordinate system, with major axes aligned with the coordinate axes, and surface $S$ given by
\begin{equation}
\frac{x^2}{a^2}+\frac{y^2}{b^2}+\frac{z^2}{c^2}=1~.
\label{E:ellipse}
\end{equation}
This is referred to as a \textit{standard ellipsoid}. The surface is parametrized by latitude-longitude coordinates $\alf\in[-\pi,\pi]$ and $\beta\in[-\pi/2,\pi/2]$, and is discretized by $n\times m$ points $\xb(\alpha_j,\beta_k)$, where $\alf_j=-\pi + j\Delta\alf,\bet_k=-\pi/2 + k\Delta\beta$, $j=0,\dots,n$, $k=0,\dots,m$, and $\Delta \alf=2\pi/n$, $\Delta\beta=\pi/m$. However, the trapezoidal rule loses accuracy for target points near the poles of such a latitude-longitude grid.  We thus use two different grids, one with the poles on the $z$-axis, a second with poles on the $x$-axis. For a given target point $\xo$, we choose the grid whose poles are furthest from $\xo$. The two grids are given by  
\begin{subequations}
\begin{align}
\xb&=\langle a\cos\alf\cos\bet,~b\sin\alf\cos\bet,~c\sin\bet\rangle\qquad \hskip 28pt \hbox{(Grid 1)}\\
\xb&=\langle a\sin\bet,~\hskip 28pt b\cos\alf\cos\bet,~c\sin\alf\cos\bet\rangle\qquad \hbox{(Grid 2)}
\end{align}
\label{E:grids}%
\end{subequations}
Grid 1 has poles at $(0,0,\pm c)$, and is discretized by $n=n_1, m=m_1$ points.
Grid 2 has poles at $(\pm a,0,0)$, and is discretized by $n=n_2, m=m_2$ points.
They are shown in figure  \ref{F:twogrids} for a sample ellipsoid, sample target point $\xo$, and sample values of $n_1,m_1$, $n_2,m_2$.

Computing the SLP and DLP using (\ref{E:changeofvars}) requires $J$ (for SLP) and $\nb J$ (for DLP). For both parametrizations (\ref{E:grids}a,b), they are given by
\begin{subequations}
\begin{eqnarray}
\nb J&=& \xb_{\alf}\times\xb_{\bet}
=abc\cos\bet\Big\langle \frac{x}{a^2},\frac{y}{b^2},\frac{z}{c^2}\Big\rangle \\
J&=&|\xb_{\alf}\times\xb_{\bet}| =
abc \cos\bet\sqrt{\frac{x^2}{a^4}+\frac{y^2}{b^4}+\frac{z^2}{c^4}}\,.\end{eqnarray}
\label{E:ellipparam}
\end{subequations}

\paragraph{Standardization of general ellipsoid}  In general, and in all our multi-ellipsoid examples, each ellipsoid $E_k$, $k=1,\dots,N$, is not in standard form, but is given as a rotation and translation of a standard ellipsoid $\tilde{E_k}$ with 
$${E}_k = \mathbf{R}_k\tilde{{E}}_k+\mathbf{s}_k,\qquad k=1,\dots,N$$
where  $\mathbf{s}_k$ is a translation vector, and $\mathbf{R}_k:=\mathbf{BCD}$ is a rotation matrix, with
$$
\mathbf{B}=\begin{bmatrix}
\cos\phi_k & -\sin\phi_k & 0\\
\sin\phi_k & \cos\phi_k & 0\\
0~~~ & 0~~~ & 1
\end{bmatrix},~~
\mathbf{C}=\begin{bmatrix}
1 & 0~~~ & 0~~~\\
0 & \cos\theta_k & -\sin\theta_k\\
0 & \sin\theta_k & \cos\theta_k
\end{bmatrix},~~
\mathbf{D}=\begin{bmatrix}
\cos\psi_k & -\sin\psi_k & 0\\
\sin\psi_k & \cos\psi_k & 0\\
0~~~ & 0~~~ & 1
\end{bmatrix}.
$$
We use rotation and translation to reduce the computation of the SLP and DLP over non-standard ellipsoids to that over the corresponding standard ellipsoid: letting 
$$
\xb=\Rb_k\xbtil+\sbk\,,\quad
\xo=\Rb_k\xotil+\sbk\,,\quad
\fb_k=\Rb_k\fbtil_k
$$
we have that
\begin{equation}
\begin{split}
{\cal{S}}_k[\fb_k](\xo)
&
=\frac{1}{8\pi\mu}\int_{S_k} \Big(\frac{\fb_k}{\rho}+\frac{(\fb_k\cdot \xbh)\,\xbh }{\rho^3}\Big) \,dS(\xb)~,\quad\\[2ex]
&=\frac{1}{8\pi\mu}\int_{\tilde S_k} \Big(\frac{\Rb_k\fbtil_k}{\rhotil}+\frac{(\Rb_k\fbtil_k\cdot \Rb_k\xbhtil)\,\Rb_k\xbhtil }{\rhotil^3}\Big) \,dS(\xbtil)~,\quad\\[2ex]
&=\frac{1}{8\pi\mu}\Rb_k\int_{\tilde S_k} \Big(\frac{\fbtil_k}{\rhotil}+\frac{(\fbtil_k\cdot\xbhtil)\,\xbhtil }{\rhotil^3}\Big) \,dS(\xbtil)~,\quad\\[2ex]
&=\Rb_k {\cal{S}}_{\ktil}[\fbtil_k](\xotil)~,
\end{split}
\end{equation}
where $
\xbhtil =\xbtil-\xotil$,
$\rhotil=|\xbtil-\xotil|$,
$\tilde{S_k}$ is the rotated surface,
and we used that $\Rb_k\xbhtil=\xbh$, $\rhotil=\rho$, $dS(\xb)=dS(\xbtil)$, and $\Rb_k\fbtil_k\cdot \Rb_k\xbhtil = \fbtil_k\cdot \xbhtil$
since $\Rb_k$ is orthogonal.
Similarly,
${\cal{D}}_k[\fb_k](\xo)= \Rb_k {\cal{D}}_{\ktil}[\fbtil_k](\xotil)$.

\subsection{The approximation $H$}
\begin{figure}
 \centering
\includegraphics[trim=0 0 0 0, clip, width=0.57\textwidth]{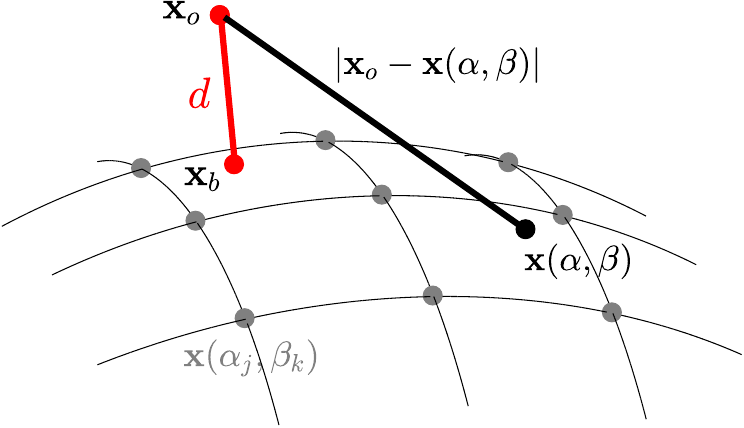}
\caption{
Sketch showing surface $\xb(\alf,\bet)$, discretization by $x(\alf_j,\bet_k)$, target point $\xo$, and corresponding distance $d$ and base point $\xb_b$.}
\label{F:sketch}
\end{figure}

We wish to approximate the integrand in 
(\ref{E:changeofvars}), 
\begin{equation}
G(\alf,\bet)= \frac{F(\xb(\alf,\bet),\xo)}{\rho^r}\omega(\xb(\alf,\bet))J(\alf,\bet)~,
\end{equation}
where $F,\omega,J$ are smooth and $\rho = |\xb(\alf,\beta)-\xo|$,
by a function $H(\alf,\bet)$ that captures the near-singularity and can be integrated exactly in closed form. 
The approximation $H$ is obtained by expanding $G(\alf,\bet)$ about a base point $(\alf_b,\bet_b)$, defined by $\xb_b=\xb(\alf_b,\bet_b)$, where $\xb_b$ is the orthogonal projection of $\xo$ onto $S$, as illustrated in figure \ref{F:sketch}. 

We begin by expanding $\xb(\alf,\beta)$ about the base point, 
\begin{equation}
\xb=\xp +\xb_{\alf}\alfh 
        +\xb_{\bet}\betah 
        +\frac{\xb_{\alf\alf}}{2}\alfh^2 
        +\xb_{\alf\bet}\alfh \betah
        +\frac{\xb_{\bet\bet}}{2}\betah^2 
        +\frac{\xb_{\alf\alf\alf}}{6}\alfh^3
        +\frac{\xb_{\alf\alf\bet}}{2}\alfh^2\betah
        +\dots
\end{equation}
where $\alfh=\alf-\alfp$, $\betah=\bet-\betp$.
In view of the orthogonality properties $(\xp-\xo)\perp \xb_{\alf}$ and $(\xp-\xo)\perp\xb_{\bet}$, we obtain the Taylor expansion of $\rho^2$, 
\begin{equation}
\rho^2=|\xb-\xo|^2
=d^2 +\calf^2\alfh^2 +2c_{\alf\bet}\alfh\betah +\cbet^2\betah^2 + \eta(\alfh,\betah)
=\rhoo^2\left[1+ \frac{ \eta}{\rhoo^2} \right]
\end{equation}
where
\begin{equation}
\rhoo^2= d^2 +\calf^2\alfh^2 + 2c_{\alf\bet}\alfh\betah+\cbet^2\betah^2,
\label{E:rhoo}
\end{equation}
and
\begin{eqnarray}
d&=&|\xb_b-\xo|,\nonumber\\ 
\calf^2&=& (\xb_b-\xo)\cdot\xb_{\alf\alf}+\xalf\cdot\xalf,\nonumber\\
c_{\alf\bet}&=& (\xb_b-\xo)\cdot\xb_{\alf\bet}+\xalf\cdot\xbet,\nonumber\\
\cbet^2&=& (\xb_b-\xo)\cdot\xb_{\bet\bet}+\xbet\cdot\xbet,
\label{E:calf}\end{eqnarray}
with all derivatives evaluated at the base point, and where
$\eta= \sum_{j+k=3 }^\infty e_{jk}\alfh^j\betah^k $ 
is the degree 3 and higher terms in the expansion of $\rho^2$.
We then expand $F/\rho^r$ with respect to $\eta/\rho_o^2$ as 
\begin{equation}
\frac{F}{\rho^r}=\frac{F}{\rho_o^r}\Big[1-\frac{r}{2}\frac{\eta}{\rho_o^2}+\frac{r(r+2)}{8}\frac{\eta^2}{\rho_o^4}- \frac{r(r+2)(r+4)}{48}\frac{\eta^3}{\rho_o^6}+\dots\Big]
\label{E:binomExpansionG}%
\end{equation}
and $F=\sum_{j+k=0}^{\infty} \tilde{c}_{jk}\alfh^j\betah^k$. 
With additional multiplication by Taylor series approximations of 
$\omega(\alf,\bet)J(\alf,\bet)$ and odd values of $r$ we obtain the expansion of $G$ in the form
\begin{subequations}
\begin{equation}
G(\alf,\bet)\approx H(\alf,\bet)=\sum c_{pqk}H_{pqk}(\alfh,\betah)~,
\end{equation}
where
\begin{equation}
H_{pqk}(\alf,\bet)=\frac{\alf^p\bet^q}{\rho_o^{2k+1}}~.
\end{equation}
\label{E:expansionH}%
\end{subequations}
$H$ consists of a truncation of the expansion (\ref{E:binomExpansionG}).
The number of terms that are kept in the expansion of $H$ in 
(\ref{E:expansionH}) are determined using the known bound for the trapezoidal rule error due to the near-singular behavior of each basis functions \cite{nitsche2021tcfd},
\begin{equation} E[H_{pqk}] =O(h^2 d^{p+q-(2k+1)}) ~,
\label{E:errorbound}\end{equation}
{where $h=\max(\Delta\alf,\Delta\beta)$.}
{We apply these error bounds for $d=O(h)$ and}
 keep all terms in $E[H]$ that are bigger than $O(h^4)$. 
 That is, we keep all terms with $2+p+q-(2k+1)\le 3$, or $p+q\le (2k+1)+1$.

The size of $F$ affects the number of terms needed in the expansion, and the accuracy needed for the density $\omega$. To illustrate, note that for the Stokes SLP, the term with $r=1$ has corresponding $F=O(1)$, while the term with $r=3$ has corresponding $F=O(d^2)$. For the DLP, the term with $r=5$ has $F=O(d^3)$. 
So the integrand has the form
\begin{equation}
G=\frac{F\omega}{\rho_o^r}(1+\frac{O(\eta)}{\rho_o^2}+\frac{O(\eta^2)}{\rho_o^4}+\dots)
\label{E:formofG}
\end{equation} where $F=O(d^m)$, some $m$, and $\eta = O(d^3)$. The $n$th term satisfies 
\begin{equation}E\Big[\frac{F\omega\,\eta^n}{\rho_o^{r+2n}}\Big]=O(h^2d^{m+3n-(r+2n)})=O(h^2d^{m-r+n}).
\label{E:errorterm}
\end{equation}
So for 4th order, we need all $n\ge 0$ up to $2+(m-r+n)< 4$, or 
\begin{equation} 
0\le n\le r+1-m.
\label{E:n_bounds}
\end{equation}
In terms of the needed basis functions,
note that for each $n$, {the term} $\displaystyle \frac{F\omega\,\eta^n}{\rho_o^{r+2n}}$ contains
{functions} $\displaystyle H_{pqk}=\frac{\alpha^p\beta^q}{\rho_o^{2k+1}}$ 
{for $p+q\geq 3n$ only, (since $\eta=O(\alfh^3)$) where
 $2k+1=r+2n$.}
For 4th order, by (\ref{E:errorbound}), we need $H_{pqk}$ for all combinations of $p, q$ up to $p+q - (2k+1) + 2 < 4$. {Thus, we require all basis functions $H_{pqk}$ with} 
\begin{equation}
3n\leq p+q\leq (2k+1)+1\equiv 2n+r+1.
\label{E:pq_bounds}
\end{equation}

{
The necessary values of $p+q$ for each $H_{pqk}$ can be obtained using (\ref{E:n_bounds}--\ref{E:pq_bounds}) and are listed in Table \ref{T:basis_count}.
For each of the relevant values of $r,m$ it lists the values of $n$ (using \ref{E:n_bounds}), and the corresponding values of $p+q$ (using \ref{E:pq_bounds}) needed to get an overall $O(h^4)$ method.}
\begin{table}[H]
\centering
{
\begin{tabular}{|c|c|c|c|c|}
\hline
~$r$~ &~$m$~ &~$n=0:(r+1-m)$~ 
&~$p+q=3n:(2n+r+1)$~&~$2k+1=r+2n $~\\ 
\hline
1&0&0& $0:2$ & 1 \\\
 & &1& $3:4$ & 3 \\
 & &2& $6:6$ & 5 \\
\hline
3&2&0& $0:4$ & 3 \\
 & &1& $3:6$ & 5 \\
 & &2& $6:8$ & 7 \\
\hline
5&3&0& $0:6$  & 5  \\
 & &1& $3:8$  & 7  \\
 & &2& $6:10$ & 9  \\
 & &3& $9:12$ & 11 \\
\hline
\end{tabular}
}
\caption{Combinations of values of $p,q,k$ for the required basis functions $H_{pqk}$ in an $O(h^4)$ correction.}
\label{T:basis_count}
\end{table}
\noindent
Thus, computing both the SLP ($r=1, m=0$ and $r=3,m=2$) and the DLP ($r=5,m=3$) to 4th order requires all basis functions with $2k+1=1,3,5,7,9,11$ with corresponding $p+q=$0:2, 0:4, 0:6, 3:8,  6:10, 9:12, which totals 179 functions.

We can use \eqref{E:errorterm} to deduce the order of the Taylor expansion required for $\omega$ and $\xb$ at the base point
for an $O(h^4)$ method. 
Let $l$ be the degree of the Taylor polynomial used to approximate $\omega$ at the base point, so that the resulting approximation  error is of order $l+1$, that is, $O(h^{l+1})$.
By \eqref{E:errorterm}, we need to include all terms in \eqref{E:formofG} satisfying $2+m+l+3n-(r+2n) = 3$, or $l= 1+(r-m)-n$. 
The largest required value of $l$ occurs when $n=0$.
Thus, for an $O(h^4)$ method, computing the SLP ($r-m=1$) requires 
a Taylor polynomial of degree $l=2$ for $\omega$, while computing the DLP ($r-m=2$) requires a Taylor polynomial of degree $l=3$.
This in turn requires $\omega$ at the base point to be accurate to fourth order, its first derivatives to third order, second derivatives to second order, and third derivatives to first order, resulting in a cubic polynomial approximation of $\omega$ at the base point with an overall error of $O(h^4)$.

Similarly, by inspecting the terms $F\omega\eta^n/\rho_0^{r+2n}$, $n\ge1$, we find that computing the SLP and DLP to $O(h^4)$, we need a polynomial approximation of $\xb$ to degree
four and five, respectively.\footnote{{For the ellipsoids considered in this paper, $\xb$ and all its derivatives are known exactly. If the surface $\xb$ were to evolve in time, it and all its derivatives at the base point could be computed locally using a 5th order polynomial interpolant to yield a 4th order method.}} 
Similar to the analysis above, let $l$ be the Taylor polynomial degree needed for $\xb$, then the largest required value of $l$ occurs when $n=1$ in \eqref{E:errorterm}, thus for the SLP and DLP we need $l=4$ and $l=5$, respectively.

From (\ref{E:formofG}--\ref{E:errorterm})  one can also deduce the size of the errors in the uncorrected trapezoidal rule. They are given by the term in the expansion (\ref{E:formofG}) with largest error, which corresponds to $n=0$. It follows that the uncorrected trapezoidal rule is $O(h^2/d)$ for the SLP and $O(h^2/d^2)$ for the DLP.

Finally, we note that the expansion (\ref{E:binomExpansionG}) requires that $\rho_0^2>0$ in (\ref{E:rhoo}) for any $d>0$, and all $\alf,\bet$, which in turn requires that the quadratic form $\displaystyle \calf^2\alfh^2 + 2c_{\alf\bet}\alfh\betah+\cbet^2\betah^2$ be positive definite (or equivalently, $\calf^2>0$, $\cbet^2>0$ and $0\le c^2_{\alf\bet}<\calf^2\cbet^2$). This holds true provided either that the target point $\xo$ is outside the ellipsoid (which is true in all examples of this paper), or that $\xo$ is inside the ellipsoid within the radius  of the osculating sphere along any tangential direction of negative curvature at the base point, see \ref{app_B}.

\subsection{Recursive evaluation of $\displaystyle \int H_{pqk}$} 
To compute the corrections $E_W[H]$, we need to accurately and efficiently compute the integrals $\displaystyle \int H_{pqk}$ over a square subdomain $W$ of the original domain approximately centered at $(\alfb,\betb)$.
With proper substitution, these integrals reduce to the form
\begin{equation}
\int_c^d\int_a^b \frac{u^p v^q}{\rho_u^{2k+1}}\,du\,dv~,\quad \rho^2_u=1+u^2+2Cuv+v^2~,
\end{equation}
where the domain is rectangular, approximately centered at the origin ($a\approx -b$, $d\approx -c$), and $\displaystyle |C|=|\frac{c_{\alf\bet}}{\calf\cbet}|<1$.  
For $C=0$, analytical expressions for all antiderivatives can be found with Mathematica. However, for $C\ne 0$, such expressions are not readily available. Instead, we obtain all antiderivatives using the recursive algorithm presented next. The values of the definite integrals are computed from the antiderivatives at each step.  The recursion, whose derivation will be presented in detail in an upcoming work, is similar to that given in \cite{helsing2013higher} for the singular integrals obtained when $d=0$.  

Define the antiderivatives 
$$ I_{pqk}= \int_c^d\int_a^b \frac{u^p v^q}{\rho_u^{2k+1}}\,du\,dv~,$$ 
for $|C|<1$ and the required values $0\le p,q\le 12$, and $0\le k\le 5$.   
In a preparatory step, evaluate all one-dimensional integrals
\begin{equation}
F_{pk}=\int_a^b \frac{u^p}{\rho_u^{2k+1}}\,du~,\quad 
G_{qk}=\int_c^d \frac{v^q}{\rho_u^{2k+1}}\,dv~,
\end{equation}
using the recursion given in \ref{app_A}, starting with known analytic expressions for $F_{00}$ and $G_{00}$.  The values of $I_{pqk}$ are then obtained as follows:

\bigskip
\noindent\textbf{Step 0:} Initialization. The definite integral of $I_{005}$ is evaluated numerically. Here we use its exact value integrated over all of $\mathbb{R}^2$, and approximate the difference between a sufficiently large domain and the actual rectangular domain using Gauss integration.
{This approach yields accurate values for the near-singular integral at low numerical cost.}
\footnote{{The sufficiently large domain used is $D=[-100,100]^2$ for which the
difference $\int_{\mathbb{R}^2\setminus D}(1/\rho_u^{11}) du\,dv <10^{-10}$.
An alternative would be to compute $I_{00k}$ for a much larger $k$, for which the difference between integrating $1/\rho_u^{2k+1}$ on $\mathbb{R}^2$ and on $[a,b]\times[c,d]$ is sufficiently small, and then use the recursion \eqref{eq:recur1_derived} to obtain $I_{005}$.}}

\medskip
\noindent\textbf{Step 1:} Evaluate $I_{00k}$, all $k<5$, using 
\begin{equation}
I_{00k}= \frac{(2k+1)I_{00(k+1)} -uG_{0k}-vF_{0k}}{2k-1}, ~k=4:-1:0.\label{eq:recur1_derived}
\end{equation}
These are computed backwards in $k$ since the integrals, of magnitude $1/R^{2k-1}$, $R\approx \max(|a|,|c|)$, increases as $k$ decreases, yielding a stable method. The reverse recursion, starting with $I_{000}$, is numerically unstable due to subtraction of like numbers.

\medskip
\noindent\textbf{Step 2:} Evaluate $I_{pqk}$ for all cases with $2k-p-q-1=0$.
Here we first use (\ref{E:recur3_derived}) 
for $p=0$, $q=2k-p-1$, all $k$,  followed by (\ref{E:recur2_derived}) 
for $q\ge 0$, $p=2k-q-1$, all $k$.
\begin{subequations}
\begin{equation}
I_{pq(k+1)}=\frac{(n-1)I_{p(q-2)k}-CpI_{(p-1)(q-1)k}+Cu^pG_{(q-1)k}-v^{q-1}F_{pk}}{(2k+1)(1-C^2)}, \label{E:recur2_derived}
\end{equation}
\begin{equation}
I_{pq(k+1)}=\frac{(p-1)I_{(p-2)qk}-CqI_{(p-1)(q-1)k}+Cv^qF_{(p-1)k}-u^{p-1}G_{qk}}{(2k+1)(1-C^2)}, 
\label{E:recur3_derived}
\end{equation}
\end{subequations}
The terms with negative indices, for example in (\ref{E:recur3_derived}) when $p=1$ or $q=0$, are not present since they are multiplied by 0. We note that all values of $I_{pqk}$ on both sides of the recursion formulas 
are of approximately equal magnitude $\sim 1/R^2$, so there is no loss of accuracy due to subtraction of like numbers.

\medskip
\noindent\textbf{Step 3:} Evaluate $I_{pqk}$ for all cases with $2k-p-q-1\ne0$. 
Here we first use (\ref{E:recur2b_derived}) 
for $q=0$, $p\ge 1$ all $k$,  followed by (\ref{E:recur2a_derived}) 
for $p\ge 0$, $q\ge 1$, all $k$.
\begin{subequations}
\begin{equation}
\begin{split}
I_{pqk} =& \big[(q-1)I_{p(q-2)k}-CpI_{(p-1)(q-1)k}+Cu^pG_{(q-1)k}-(1-C^2)u^{p+1}G_{qk}\\
& -(v^{q-1}+(1-C^2)v^{q+1})F_{pk}\big]/\big[(1-C^2)(2k-p-q-1)\big]
\end{split}
\label{E:recur2a_derived}
\end{equation}
\begin{equation}
\begin{split}
I_{pqk} =& \big[(p-1)I_{(p-2)qk}-CqI_{(p-1)(q-1)k}+Cv^qF_{(p-1)k}-(1-C^2)v^{q+1}F_{pk}\\
&-(u^{p-1}+(1-C^2)u^{p+1})G_{qk}\big]/\big[(1-C^2)(2k-p-q-1)\big]
\end{split}
\label{E:recur2b_derived}
\end{equation}
\end{subequations}
where, as above, all $I_{pqk}$ with negative indices are set to zero.

\subsection{Trapezoidal Rules $T^4$ and $T^6$}
\label{S:trapez}
The trapezoidal rules of fourth and sixth order used here to compute double integrals over rectangular domains discretized by a uniform mesh are obtained by applying the Euler-Maclaurin expansion for 1D integrals in each direction. We have that
\begin{align}
\int_c^d\int_a^b f(\alf,\bet)\,d\alf\,d\bet \approx &
\delbet\delalf \sumtrap{k}{m} \sumtrap{j}{n} f(\alf_j,\bet_k) \\
&- \frac{\delbet\delalf^2}{12} \sumtrap{k}{m} \frac{\partial f}{\partial \alf}(\alf,\bet_k) {\bigg|}_{\alf=a}^{\alf=b}
- \frac{\delbet^2\delalf }{12} \sumtrap{j}{n} \frac{\partial f}{\partial \bet} (\alf_j,\bet) {\bigg|}_{\bet=c}^{\bet=d}
\nonumber\\
&
{+ \frac{\delbet\delalf^4}{720} \sumtrap{k}{m} \frac{\partial^3 f}{\partial \alf^3}(\alf,\bet_k) {\bigg|}_{\alf=a}^{\alf=b}}
{+ \frac{\delbet^4\delalf}{720}
\sumtrap{j}{n} \frac{\partial^3 f}{\partial \bet^3} (x_j,\bet) {\bigg|}_{\bet=c}^{\bet=d} }\nonumber\\
&
{+ \frac{\delbet^2\delalf^2 }{12^2} \frac{\partial^2 f}{\partial \alf\partial \bet}(\alf,\bet) {\bigg|}_{\alf=a}^{\alf=b} {\bigg|}_{\bet=c}^{\bet=d}}
\nonumber
\end{align}
where the prime denotes that the first and last term in the sum are weighted by 1/2. 
The first term on the right hand side is of order 1. The second and third are of order $h^2$, where $h=\max{\delalf,\delbet}$.  The last three are all the terms in the expansion of order $h^4$.
The 2nd, 4th, and 6th order trapezoidal rules $T_D^2[f], T_D^4[f], T_D^6[f]$ are thus given by the first one, three and six terms on the right hand side, respectively. As described next, 
the method presented in this paper uses $T^4$ and $T^6$.
For ease of notation, function values at the $(j,k)$th grid point will be denoted below by indices $j,k$, for example, $f(\alf_j,\bet_k)=f_{j,k}$.

\subsection{On the size $n_w$ of the window $S$}
\label{S:window}
The method consists of correcting the trapezoidal approximation of
near-singular integrals 
by integrating a local approximation exactly over a subdomain $W$ centered on the near-singularity.
Specifically, the approximation (3.2, 3.3), repeated here for clarity, is
\begin{equation}
\begin{split}
\int_D G& = \int_{D\backslash W}G +\int_{W}(G-H)+\int_W H\\
&\approx T_{D\backslash  W}[G] + T_{W}[G-H] +\int_W H\\
&=T_{D}[G] + E_{W}[H]
\label{E:methodagain}
\end{split}
\end{equation}
Note that the term $T_{W}[G-H]$ requires sampling both $G$ and $H$  at the same grid points.  The subdomain $W\subset D$ is thus chosen to be an $2n_w\times 2n_w$ subgrid of the original grid, roughly centered on $(\alfb,\betb)$, $W=[\alf_{j_0-n_w},\alf_{j_0+n_w}]\times[\bet_{k_0-n_w},\bet_{k_0+n_w}]$ such that $(\alf_{j_0},\bet_{k_0})$ is the nearest grid point to $(\alfb,\betb)$. Since $W$ is rectangular, $\int_{W} H$ can easily be evaluated using the recursive formulas given above.

Next we describe how to choose the size of the window $W$.
The term $E_W[H]$ captures the error of the trapezoidal approximation of near-singular integrals. This term is bounded by (\ref{E:errorbound}). 
However, it also contains the trapezoidal error due to derivatives at the boundary $\partial W$ of $W$. If the window size $n_w h$ is small, the derivatives of $H_{pqk}$ are large on $\partial W$. The size of these errors is best captured by the error in integrating $H$ over the complement, $E_{D\backslash W}[H]$. This error is dominated by the boundary derivatives, since $H$ is regular (and not near-singular) outside $W$.
Furthermore, the error depends on the quadrature rule used. 

\begin{figure}
\begin{center}
\includegraphics[trim= 0 0 0 0, clip, width=3.8truein]{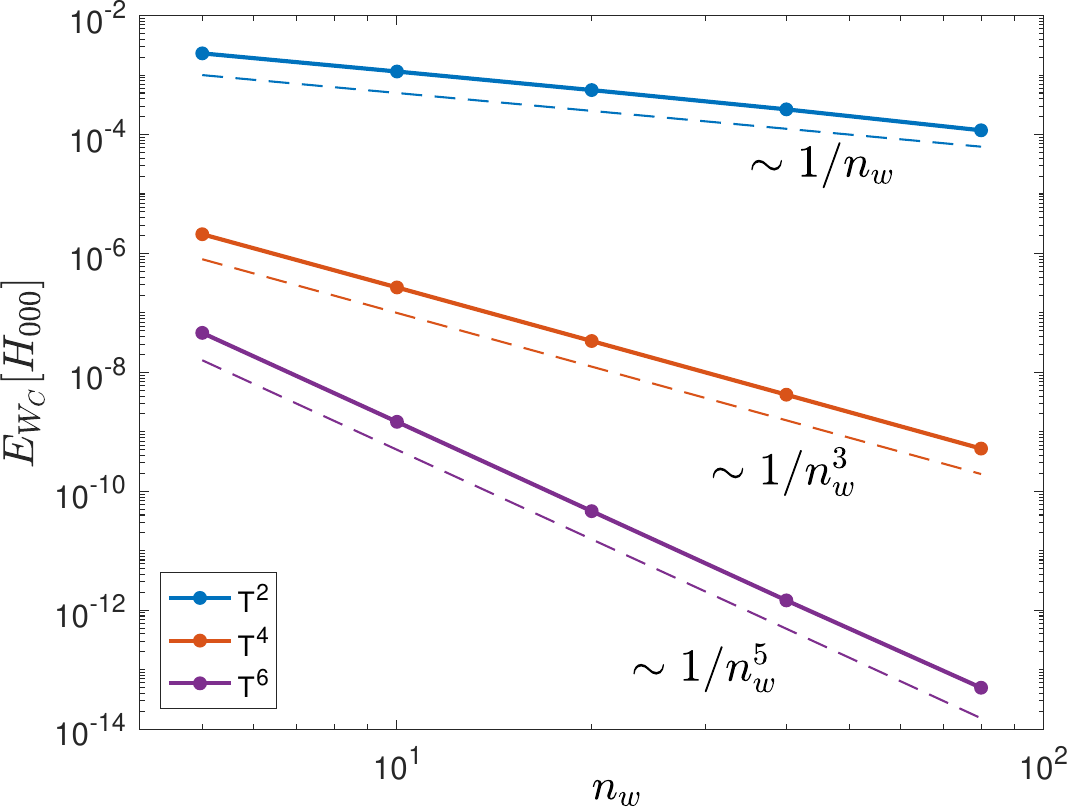}
\end{center}
\caption{{Error in the numerical integration of $H_{000}=1/\sqrt{x^2+y^2+d^2}$ with $d=0.001$ over the complement $W_C=D\setminus W$, where $W=[-\nwin h,\nwin h]\times [-\nwin h,\nwin h]$,  $D=[-10,10]\times [-10,10]$} and $h=0.025$, as a function of $n_w$, using the 2nd, 4th and 6th order Trapezoidal rules $T^2, T^4, T^6$, as indicated. }
\label{F:window}
\end{figure}

Figure \ref{F:window} plots the error in integrating a sample basis function,  
$H_{000}(x,y)=1/{\sqrt{x^2+y^2+ d^2}}$, $d=0.001$,  over the complement $D\backslash W$ of a window $W=[-\nwin h,\nwin h]\times [-\nwin h,\nwin h]$ in the domain $D=[-10,10]\times [-10,10]$, discretized by $h=0.025$, using the three quadratures defined above, $T^2$, $T^4$ and $T^6$, as a function of $n_w$. The figure shows that smaller values of $n_w$ yield larger errors, as expected, and that the errors are smaller using higher order trapezoidal rules. {Equivalently, for a given desired error size, a smaller window can be used if we use a higher trapezoidal rule. Here, we use $T^6$ to compute $E_W[H]$, that is, we use}
\begin{equation}
\int_D G \approx T^4_{D}[G] + E^6_{W}[H]
\label{E:finalapprox}
\end{equation}
where $\displaystyle E_W^6[H]=\int_W H - T_W^6[H]$. 
This amounts to evaluating (\ref{E:methodagain}) using $T^4$, and adding $[E^6_W[H]-E^4_W[H]]$. 
That is, we use
\begin{equation}
\int_D G \approx T^4_{D}[G] + E^4_{W}[H] +(E^6_W[H]-E^4_W[H])
\end{equation}
Note that the added term $(E^6_W[H]-E^4_W[H])$ is $O(h^4)$ so this addition does not affect the order of the method. It does affect the required window size, since the edge derivatives in $E^6$ are smaller than in $E^4$. Thus $E^6$ permits using a smaller window size $n_w$.

A scaling rule for  the window size can be obtained as follows: per figure 4, using $T^6$ on the window yields $O(1/\nw^5)$ errors. For these errors to be $O(h^4)$, with $h=O(1/n)$, one requires $1/\nw^5\sim 1/n^4$, or $\nw =O(n^{4/5})$. 
The additional cost to compute the correction is a sum over $(2n_w)^2$ points for 179 basis functions.
This cost is $O(n^{8/5})=O(n^{1.6})$ and grows slower than the trapezoidal rule's $O(n^2)$ cost.
The cost can be further reduced if, for example, $T^6$ is replaced by $T^8$, in which case $\nw=O(n^{4/7})$, then the correction cost is reduced to $O(n^{1.14})$.
We use $T^6$ for simplicity while allowing $\nw$ to grow slower than $n$.

Finally, note that the sample basis function used in figure 4, with numerator 1 corresponding to $p=q=0$, captures the worst case of all basis functions with $p+q\ge 0$. According to Table \ref{T:basis_count}, in all cases $p+q\le (2k+1)+1$, so the integrand grows at most linearly. Since the trapezoidal rule is exact for linear functions, the numerical integration error equals that of integrating the function minus a linear far-field. This difference can be shown to be smaller than the $p=q=0$ case, with smaller edge derivatives. Thus the numerator 1 captures the worst case scenario.

The size of the window used in the results shown below ranges from $n_w=5$ for $n\le 80$ and $n_w=9,15,26$ for $n=160,320,640$, respectively. These values correspond to $\nw \approx 0.15 n^{4/5}$. 

The focus of the present work is the presentation of the method and its accuracy, and the runtimes in all cases were very short. Thus we did not focus on efficient implemention of the method. However, our timings are consistent with the estimates given above. Using $T^6$, the expected numerical cost is $O(n^{1.6})$. For our results presented in \S 4.2 below, the times taken to compute the corrections with $n_w=5,9,15, 26$ are well approximated by $10^{-5}n^{1.52}+0.14$ seconds per 1000 target points, 
on a single core of a 2023 MacBook Pro M2 Pro. With $n=80, n_w=5$, the corresponding observed throughput was 4000 corrections per second. Reducing the numerical cost using higher order corrections, vectorization and parallelization is the subject of future work.

\subsection{Minimizing Roundoff Error}

\begin{figure}
\centering
\includegraphics[trim= 16 12 0 17, clip, width=2.8truein]{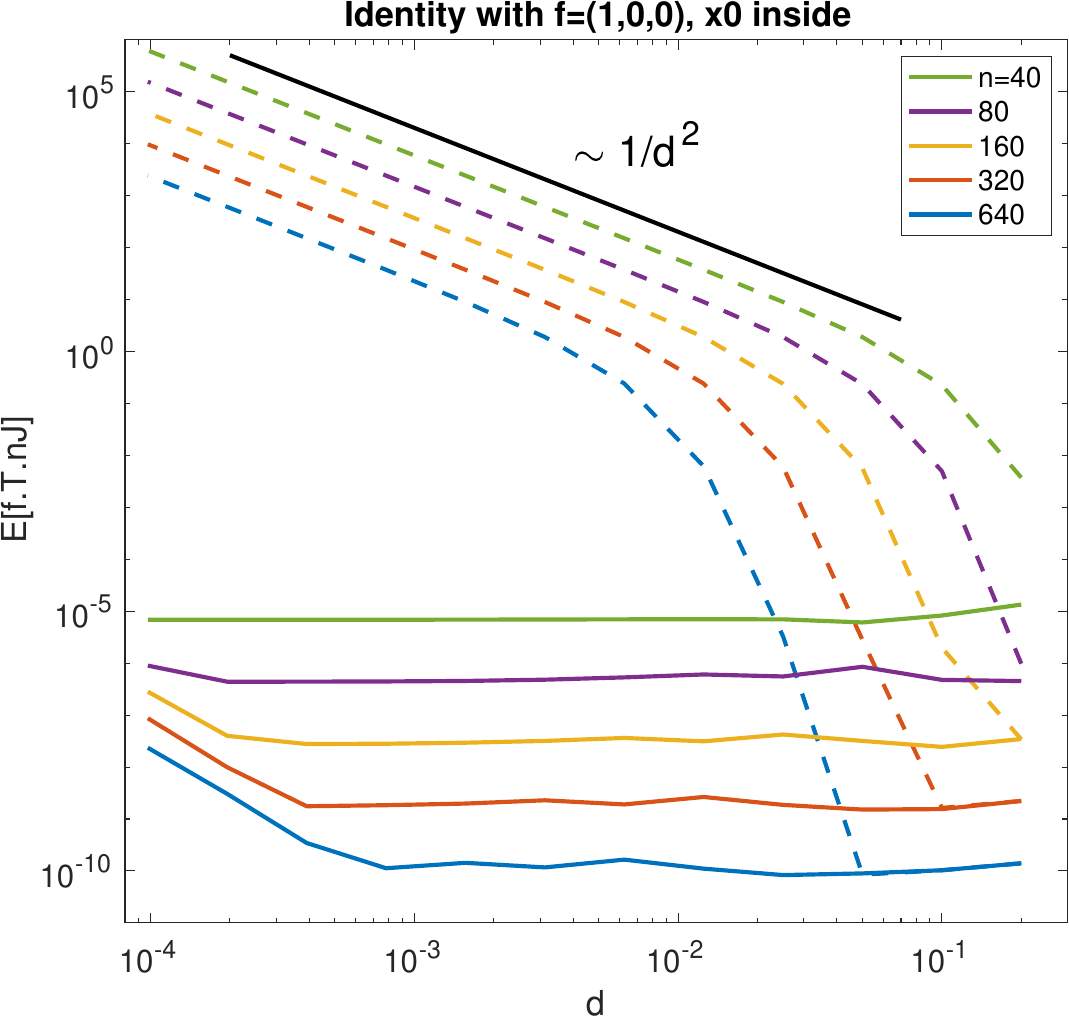}\qquad
\includegraphics[trim= 16 12 0 17, clip, width=2.8truein]{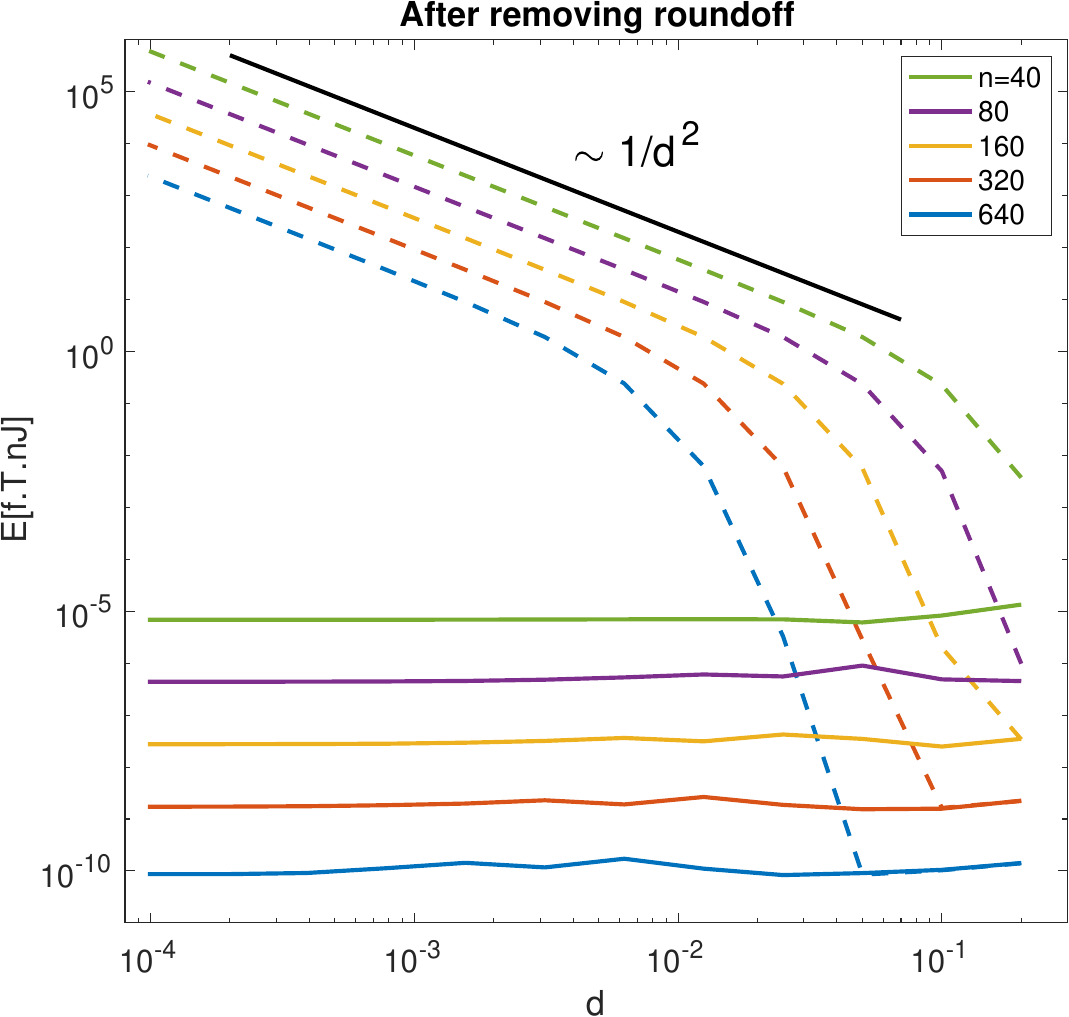}\\
\vskip-2.6truein {\footnotesize (a)}\hskip2.97truein{\footnotesize (b)}\hskip 2.9truein\hbox{ } \\
\vskip.5truein
\rotatebox{90}{\scriptsize Error in $D[f](\mathbf{x_o})$}
\hskip5.9truein\hbox{ }\\
\vskip0.8truein {\footnotesize d\hskip3.1truein d}\hskip .0truein\hbox{ } 
\caption{
Maximal error in $D[\fb](\xo)$ with $\fb=(1,0,0)$ and $\xo$ inside sphere at distance $d$ from the surface, vs $d$.  The integral is computed with the uncorrected trapezoidal rule $T^4[G]$ (dashed curves) and after adding the correction $E^6[H]$ (solid curves), with several values of $h = 2\pi/n$, $n$ as indicated.  (a) Corrected results show loss of accuracy due to subtraction of like numbers when $d \ll 1$. (b) Corrected results are $O(h^4)$, uniformly in $d$, when the punctured trapezoidal rule is used.} 
\label{F:roundoff}
\end{figure}

The method consists of obtaining accurate results by adding a correction to a bad result. If the bad result is large, the correction is large, resulting in subtraction of two large numbers to obtain a more accurate much smaller number. This leads to loss of significant digits. 
This effect is more visible in calculations of the DLP, since the maximum uncorrected errors are larger, $O(h^2/d^2)$, than for the SLP, where they are $O(h^2/d)$. Figure \ref{F:roundoff}(a) shows the error in calculating the double layer $D[\fb](\xo)$ for a sphere with $\fb=(1,0,0)$ and $\xo$ inside the sphere, for which the integral is known to be $(-1,0,0)$. It shows both the uncorrected error (dashed curves) and the error after adding the correction $E^6[H]$ (solid curves), for a range of values of $h=2\pi/n$, where $n$ is as indicated. For each $h$, it shows the maximum errors at distance $d$ from the surface. The uncorrected error behaves as expected, like $h/d^2$, for each $h$. The corrected errors are up to $10^{11}$ times smaller, and are almost constant in $d$. However, there is  loss of accuracy for $n\ge 80$ as $d\to 0$.
This loss of accuracy is caused by subtraction of like numbers in (\ref{E:finalapprox}), when the values of $T^4[G]$ and $E^6[H]$ are large and almost equal in magnitude, but of opposite sign.
Specifically, if we denote the trapezoidal rule as
\begin{equation}
T[G] = \sum_{j,k} G_{j,k}h^2 
\end{equation}
then the loss of significance is due to the subtraction $G_{j_0,k_0}h^2-H_{j_0,k_0}h^2$ at the grid point $\xb_{j_0,k_0}$ closest to the target point $\xo$.

This issue can be resolved by omitting the terms $G_{j_0,k_0}h^2$ and $H_{j_0,k_0}h^2$ from the trapezoidal rule summations.\footnote{For the DLP, an alternative approach to reducing the loss of significance is by a global singularity subtraction using the identity \cite[Eq.6.4.5]{pozrikidis1992boundary} on a closed surface. However, our proposed resolution is simpler, requires only a local modification, and applies to any integral of the form \eqref{E:genform} and even on open surfaces.} In practice, this amounts to replacing the trapezoidal rule in $T_D[G]$ and $T_W[H]$ by a ``punctured trapezoidal rule,'' which is a common strategy used in singular quadratures, see \cite{wu2021corrected} for example. Figure \ref{F:roundoff}(b) shows the removal of roundoff error by this process for the example introduced in figure \ref{F:roundoff}(a). The error in figure \ref{F:roundoff}(b) is $O(h^4)$, decaying from $10^{-5}$ to $10^{-10}$ as $h$ decreases by a factor of 16. In particular, the error converges \textit{uniformly} in $d$.
We follow this approach when $\alfh$, $\betah$, and $d$ are sufficiently small. The details are given below in \S 3.8, which lists the specific steps taken in the numerical implementation.

We now present the analysis to justify the use of the punctured trapezoidal rule. 
First note that when $(j,k)=(j_0,k_0)$ and $\alfh,\betah\approx 0$, we have 
$\rhoo\approx d$ in \eqref{E:rhoo}.
For a function $G$ of the form \eqref{E:formofG} with 
$F\sim d^m$ for some $m<r$, we have
\begin{equation}
    G_{j_0,k_0}h^2 \sim \frac{F\omega}{\rhoo^{r}}\Big|_{(\alpha_{j_0},\beta_{k_0})}h^2 \sim d^{m-r}h^2 \to \infty ~\text{as}~d\to0.
    \label{E:Gj0k0_asymp}
\end{equation}
By the expansion \eqref{E:expansionH} and the truncation condition \eqref{E:pq_bounds}, $G-H$ is a remainder of the expansion containing basis functions $H_{pqk}$ that satisfy $p+q-(2k+1)\geq2$, such that
\begin{equation}
    H_{pqk}\Big|_{(\alpha_{j_0},\beta_{k_0})} = \dfrac{\alfh^p\betah^q}{\rhoo^{2k+1}}\Big|_{(\alpha_{j_0},\beta_{k_0})} = O(d^{p+q-(2k+1)}) = O(d^2)
\end{equation}
Combining this fact with the Taylor Theorem one can show that $(G-H)_{j_0,k_0}=O(d^2)$,
which implies
\begin{equation}
    (G-H)_{j_0,k_0}h^2 = O(d^2h^2).
    \label{E:G-Hj0k0_asymp}
\end{equation}
Thus \eqref{E:Gj0k0_asymp} and \eqref{E:G-Hj0k0_asymp} shows that $G_{j_0,k_0}h^2-H_{j_0,k_0}h^2$ is indeed a subtraction of two large numbers that cancel to $0$ as $d\to0$, whereas ignoring these terms in $T_D[G]$ and $T_W[H]$ only introduces an  $O(d^2h^2)= O(h^4)$ error, since $d=O(h)$.

On the other hand, when $(j,k)\neq(j_0,k_0)$, evaluating $(G-H)_{j,k}h^2 = G_{j,k}h^2 - H_{j,k}h^2$ does not lead to loss of significance. Note that in this case $|(\alfh,\betah)| = |(\alf_{j}-\alf_b, \beta_{k}-\beta_b)| \ge h/2$, thus $\rhoo \geq \sqrt{d^2+c\cdot(h/2)^2}> \frac{\sqrt{c}}{2}\,h$ where $c>0$ is the minimum singular value of a positive definite matrix associated with the quadratic form in \eqref{E:rhoo}; see also \ref{app_B}. Thus the size of $G_{j,k}h^2$ (and likewise $H_{j,k}h^2$) is estimated to be
\begin{equation}
    G_{j,k}h^2 \sim \frac{F\omega}{\rhoo^{r}}\Big|_{(\alpha_{j},\beta_{k})}h^2 = O(d^mh^{-r}h^2) = O(h^{m-r+2}).
\end{equation}
So $G_{j,k}h^2$ and $H_{j,k}h^2$ are bounded for all $d$ if $m-r+2\geq 0$, which is the case for the integrals considered in this paper.

\subsection{Implementation}

The following steps are taken to implement the method described above to compute the SLP and DLP 
${\cal{S}}[\fb](\xo)$,
${\cal{D}}[\fb](\xo)$
over the surface of an ellipsoid $E$ with major axes $a,b,c$ in arbitrary location. Let $a$ denote the largest major axis.
{Let $h=\max(\Delta\alpha, \Delta\beta)$, where $\Delta\alpha, \Delta\beta$ are the chosen mesh size.}
We also include references to the earlier relevant sections in each step.
{Step 0 is precomputed and applies to all target points.
Steps 1-5 are taken for each individual target point $\xo$.}

\begin{itemize}
\item[{\textbf{Step 0.}}] {Precomputation: Standardize $E$ and set density.}
\item[{\textbf{0a.}}] {Rotate ellipse $E$ and $\xo$ to standard position $\tilde{E}$ (\S3.2).}
\item[{\textbf{0b.}}] {Set the density $\tilde{f}$ corresponding to the prescribed far field on both grids (\S 2). Since the flow is steady, the density stays constant in time.}
\\
\item[\textbf{Step 1.}] Determine grid.\\
Determine which grid to use: Check if $\xotil$ is closer to poles of grid 1 or of grid 2. Choose the grid whose poles it is further away from (\S3.2). 
\item[\textbf{Step 2.}] 
Determine if correction is needed. If not, compute $T[G]$. 
{If yes, proceed with steps 3-5}.
\item[\textbf{2a.}] 
The correction is applied if the distance $d<6\Delta s$, where $\Delta s$ is the physical grid size at the base point. Whether this criterion holds is determined in 2 steps. The choice of the factor 6 is addressed below, in \S 4.1.
\begin{itemize}
\item[(i)] 
We first find a lower bound $d_{min}$ for $d$, determined by the minimum distance between the ellipsoid $E$ and  the ellipsoid $\tilde E$ through the target point, as well as a upper bound $\Delta s_{max}$ for $\Delta s$.\footnote{{$d_{min}= (r-1)\min(a,b,c)$, where $r^2=x_0^2/a^2+y_0/b^2+z_0^2/c^2$, $\Delta s_{max}=\max(a,b,c)h$}} If $d_{min}/\Delta s_{max}>6$ then $d/\Delta s>d_{min}/\Delta s_{max}>6$ and no correction is needed. Set $\texttt{correct=false}$. Else, proceed to (ii).
\item[(ii)] Find the orthogonal projection $\xbb$: starting from a nearby grid point, the nonlinear equation $(\xb_b-\xo)=d\nb$, such that $\xb\in S$, is solved using Newton's method. Find $d$ and $\Delta s$.  If $d<6\Delta s$, set \texttt{correct=true}. Else, set \texttt{correct=false}.
\end{itemize}
\item[\textbf{2b.}] If \texttt{(not correct)} set $\int G=T^4[G]$ (\S3.5). Else, continue with steps 3-5.

\item[\textbf{Step 3.}] Remove roundoff error if necessary. (\S3.7)
\item[\textbf{3a.}] Roundoff error is large if $d/\Delta s$ is small and the base point is close to a grid point. Here, we determine whether $d<\Delta s/4$ and $\min\limits_{j,k}|\xb_b-\xb_{j,k}|<\Delta s/4$. If so, set \texttt{roundoff=true}.
\item[\textbf{3b.}] 
If \texttt{(not roundoff)}, compute the trapezoid sum $T[G]=T^4[G]$.
If \texttt{(roundoff)}, compute the punctured sum obtained by skipping $(j_0,k_0)$ term (\S3.7).
\item[\textbf{Step 4.}] Compute $E[H^6_{pqk}]$ and $c_{pqk}$
\item[\textbf{4a.}] Compute all derivatives of $\xb,F,J,\fb$ at $(\alfb,\betb)$.  Compute the coefficients $c_{\alf}^2$, $c_{\alf\bet}$, $c_{\bet}^2$.
The derivatives of $\xb, F, J$ at the base point are given analytically. The derivatives of $\fb$ are obtained by bi-cubic interpolation of the given values at the grid points.
\item[\textbf{4b.}] Compute $E[H_{pqk}]$. Note that the functions $H_{pqk}$ depend only on the base point $(\alfb,\betb)$, the distance $d$ and the coefficients
$c_{\alf}^2$, $c_{\alf\bet}$, $c_{\bet}^2$.
\begin{itemize}
\item[(i)] 
Set the window $W$ using a pre-determined value of $n_w$ (\S3.6).
\item[(ii)] Compute $\int_W H_{pqk}$ using the recursive algorithm (\S 3.4).
\item[(iii)] Compute $T^6_W[H_{pqk}]$. If \texttt{(roundoff)}, skip the $(j_0,k_0)$ term (\S 3.5).
\item[(iv)] 
Per (\ref{E:finalapprox}), set $$E^6[H_{pqk}]=T^6_W[H_{pqk}]-\int_W H_{pqk}$$
\end{itemize}
\item[\textbf{4c.}] Compute all coefficients $c_{pqk}$. These are obtained by multiplying the corresponding coefficients of the expansions of $F$, $J$, $f$ and $\eta^l$ for the proper value of $l$.
\item[\textbf{Step 5.}] Set $\int G=T[G]+E^6[H]~, \hbox{ where } E^6[H]=\sum_{pqk}c_{pqk}E^6[H_{pqk}]$ (\S 3.1).
\end{itemize}

{A preliminary Fortran implementation of the proposed near-singular correction for the Stokes SLP and DLP, as well as the source code used to produce the examples in the numerical results section, is available at: \texttt{https://github.com/monikanitsche/3dstokes-multi}.}

\section{ Numerical results} 
In this section we apply the corrected quadrature $T^4[G]+E^6[H]$ to compute Stokes flow past several objects, and compare with results using the uncorrected trapezoidal rule $T^4[G]$. 
The goal is to show the accuracy of the corrected values, and the severity of the errors in the absence of corrections.

\subsection{Stokes flow past a sphere}
\begin{figure}
 \centering
\includegraphics[trim=70 0 70 25, clip, width=0.45\textwidth]{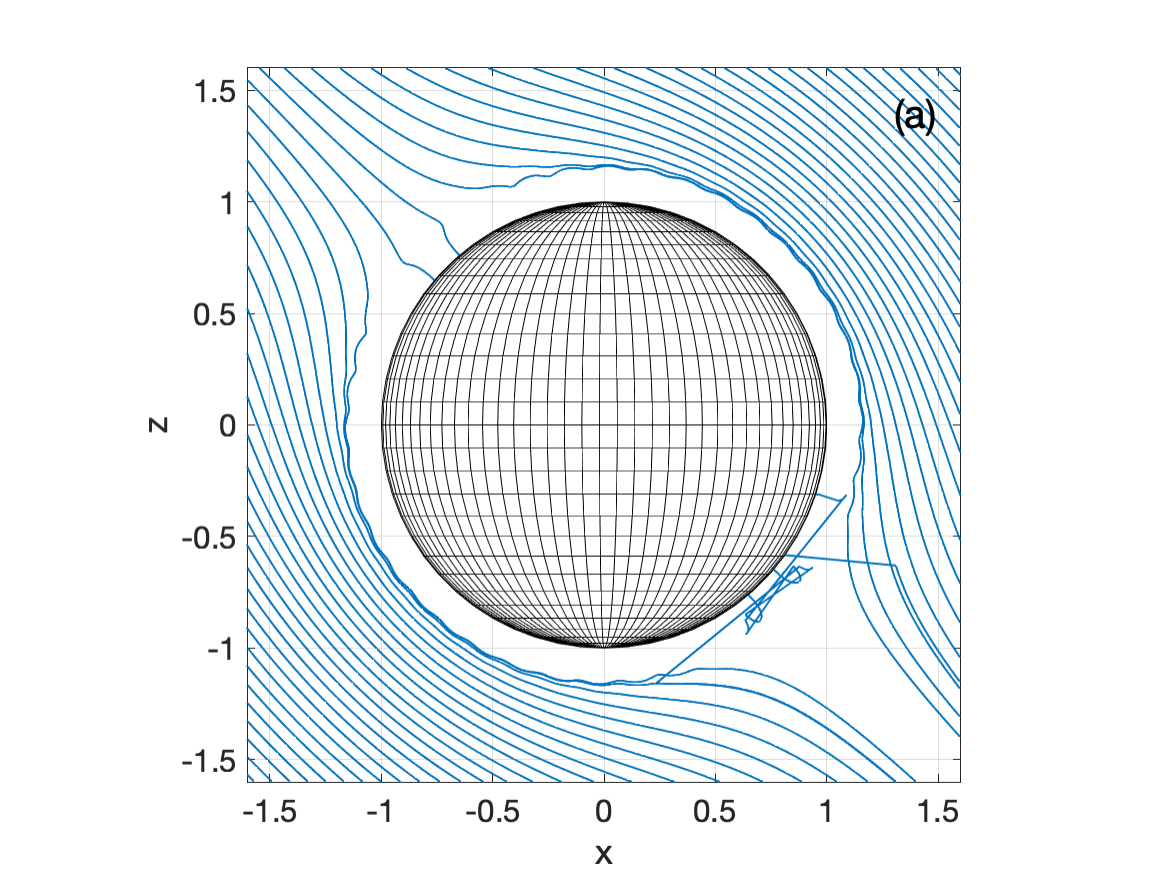}
\includegraphics[trim=70 0 70 25, clip, width=0.45\textwidth]{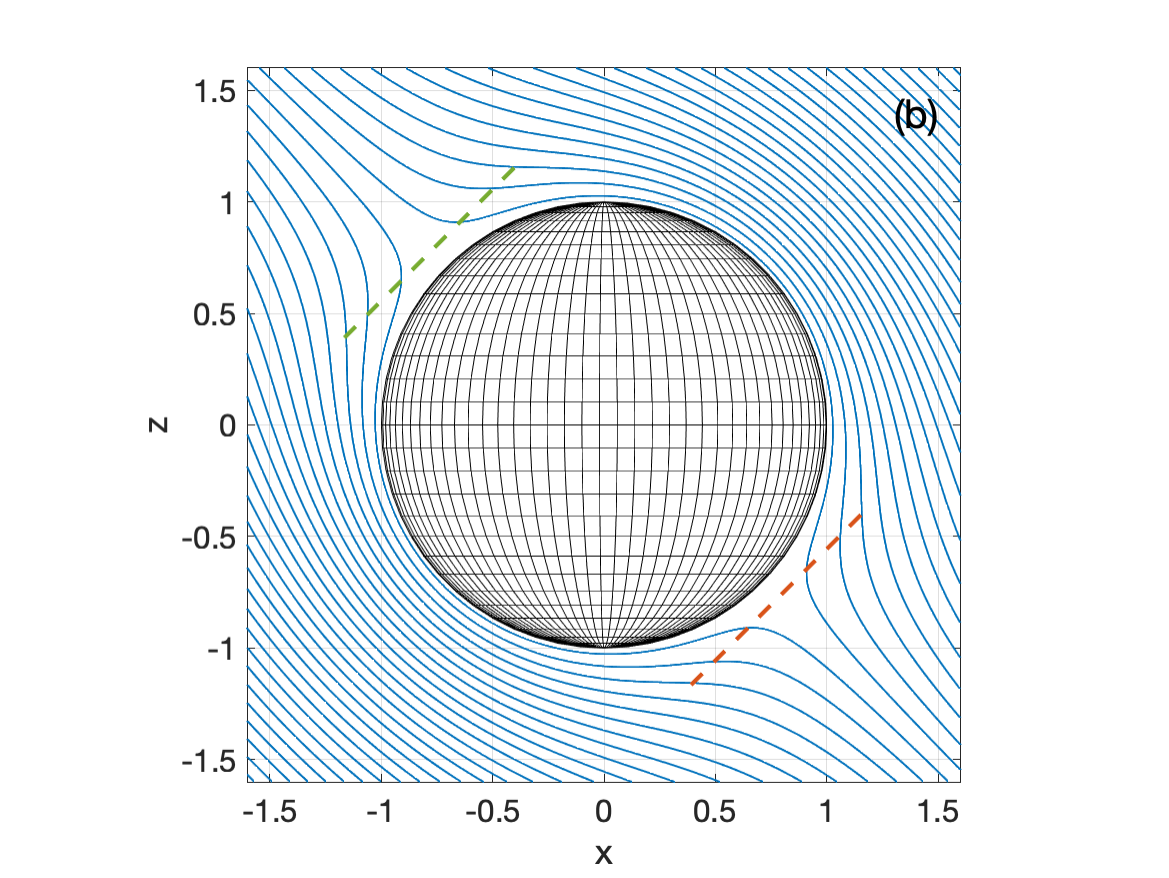}
\caption{
Streamlines in the plane $y = 0$, obtained by integrating the velocity field along particle paths with RK4 and $\Delta t = 0.1$. The velocity is computed with $n=40$, without corrections in (a) and with corrections in (b).  The green and red dashed lines in (b) are included for later reference.}
\label{F:1Sstream}
\end{figure}

Our simplest example is flow past a sphere of unit radius with constant far field velocity $\Ubinf$ of unit magnitude. Here, the density $\fb=-1.5\Ubinf$ and the fluid velocity are known explicitly \cite{ockendon1995}, making comparison with exact solutions possible. In addition, the parametrizations (\ref{E:grids}) are orthogonal, with $c_{\alf\bet}=0$ (\ref{E:calf}), thus reducing sources of error. Finally, for constant density and target points outside the sphere the DLP vanishes. We thus set
\begin{equation}
\ub(\xo)={\cal{S}}[\fb](\xo)+\Ubinf
\label{E:velo1S}
\end{equation}
and solely test the SLP.  Results solely for the DLP were shown earlier, in figure\ \ref{F:roundoff}.  Throughout this section, the mesh used is defined by $n$, where $n_1=n_2=n$, $m_1=m_2=n/2$. To avoid grid imprinting, the far field is chosen to be at a $45^o$ degree angle from the equator, $\Ubinf=(1,0,-1)/\sqrt{2}$.

Figure \ref{F:1Sstream} plots the streamlines past the sphere in the cross section $y=0$. They are computed by integrating the velocity field for a finite number of initial particle positions using the 4th order Runge Kutta method (RK4) with $\Delta t=0.1$, which is sufficiently small that results are resolved in time. The velocity along particle paths is computed using $n=40$, without corrections in panel (a), and with corrections in panel (b). 
The uncorrected velocity in (a) has large errors near the boundary, and the streamlines within that region are not resolved. Most streamlines stay outside a layer around the sphere and do not properly approach it, while some streamlines enter the sphere.
{This behaviour of the uncorrected velocity field is characteristic of all cases presented below.} 
The streamlines with the corrected velocity, in figure \ref{F:1Sstream}(b), are well resolved, even close to the boundary.

\begin{figure}
 \centering
\includegraphics[trim=0 0 0 0, clip, height=0.39\textwidth]{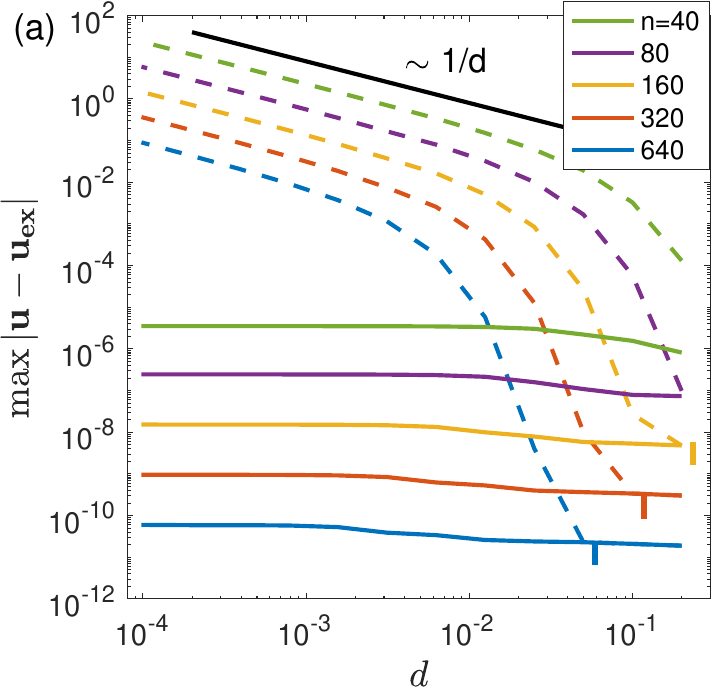}
\includegraphics[trim= 0 0 295 25, clip, height=0.39\textwidth]{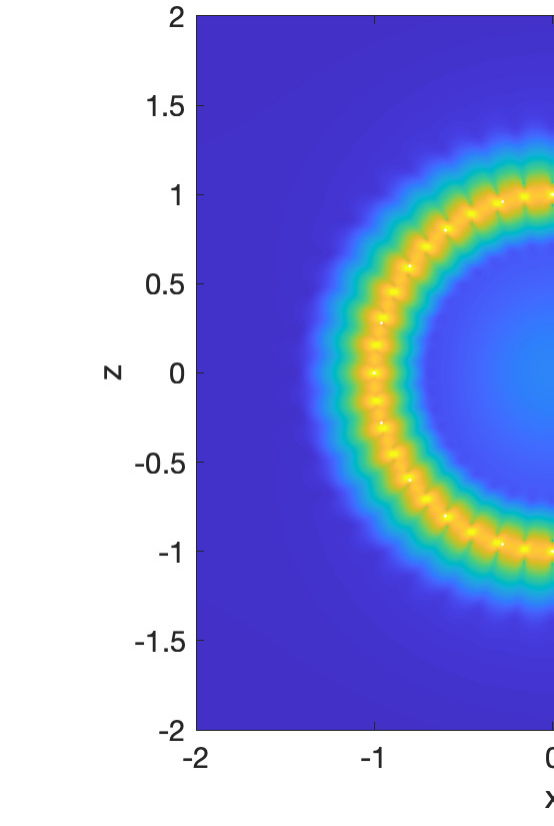}
\hskip-0.06in\includegraphics[trim= 265 0 32 25, clip, height=0.39\textwidth]{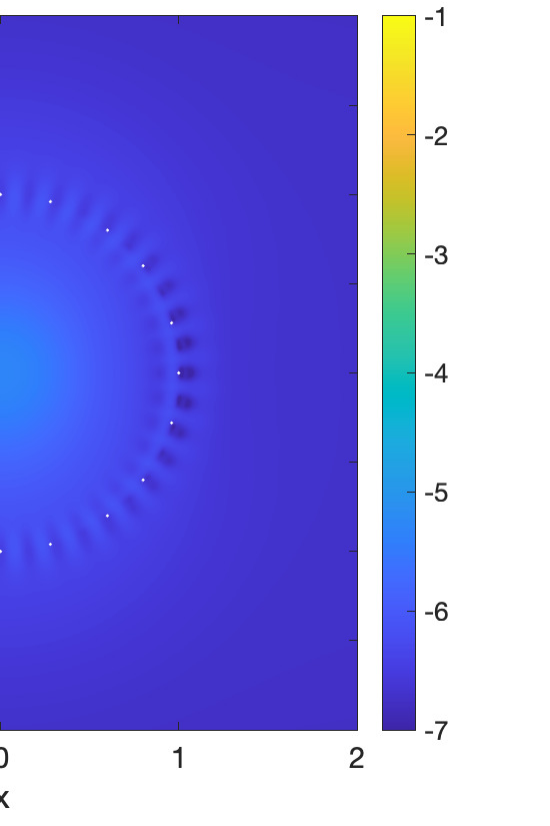}
\vskip-6.6truecm\hskip0.truecm {\footnotesize (b)}
\vskip-0.68truecm \hskip14.08truecm {\tiny $+$}
\vskip7.truecm

\caption{ (a) Maximal error in the velocity past the sphere at distance $d$ from the surface, vs $d$.  The velocity is computed without corrections  (dashed curves) and with corrections (solid curves), using the indicated values of $n$. A line proportional to $1/d$ is shown as reference.
For $n\ge 160$, a short vertical line indicates the distance $d_{max}=6\Delta s_{max}$ past which no corrections were applied.
{(b) Distribution of error in the $y=0$ plane, with $n = 40$, computed without corrections (left half) and with corrections (right half), on a logarithmic scale.}
}
\label{F:1Sconv}
\end{figure}

{Figure \ref{F:1Sconv} plots the magnitude and distribution of the  error in the computed velocity. } Figure \ref{F:1Sconv}(a) plots the maximum error over all points at a fixed distance $d$ from the surface, for a range of values of $d$ and several mesh sizes $h=2\pi/n$, with $n$ as indicated. The maximum error is evaluated by sampling the error over $320^2$ points in $(\alf,\bet)\in[0,\pi/4]\times[0,\pi/4]$, for each $d$.  The uncorrected errors (dashed curves), satisfy the expected scaling $O(h^2/d)$.  The corrected errors (solid curves) are several orders of magnitude smaller and are uniformly $O(h^4)$ as $d\to 0$. 
The maximal distance $d_{max}=6\Delta s_{max}$ (as described in Step 2 in \S3.8) above which no correction is applied is indicated by a short vertical line for $n\ge 160$. For smaller $n$, $d_{max}$ lies outside the domain of the figure. The figure shows that for $d\approx d_{max}$ the corrections no longer improve the standard fourth order trapezoid rule, making no further correction necessary for larger $d$.

Figure \ref{F:1Sconv}(b) shows the distribution of the error in the $z=0$ plane, computed with fixed $n=40$. The left side plots the uncorrected error. These are largest near the gridpoints. The right side shows the sigificant error reduction after adding the correction. The ``$+$'' sign in the legend denotes that the maximal errors are bigger than $10^{-1}$.

\begin{figure}
\centering
\includegraphics[trim= 0 0 0 0, clip, width=2.6truein]{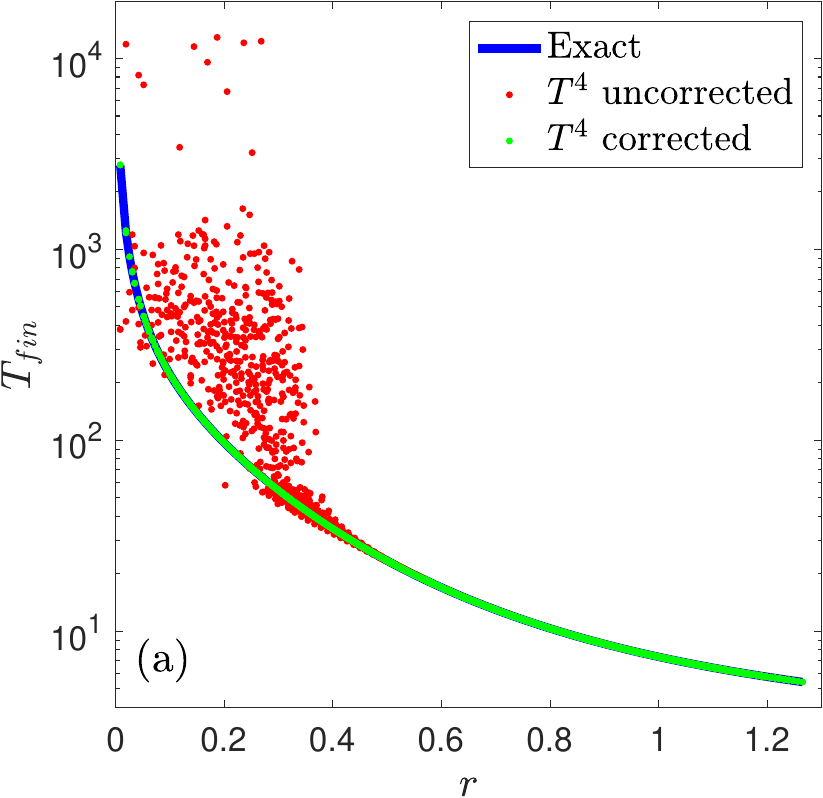}
\includegraphics[trim= 0 0 0 0, clip, width=2.6truein]{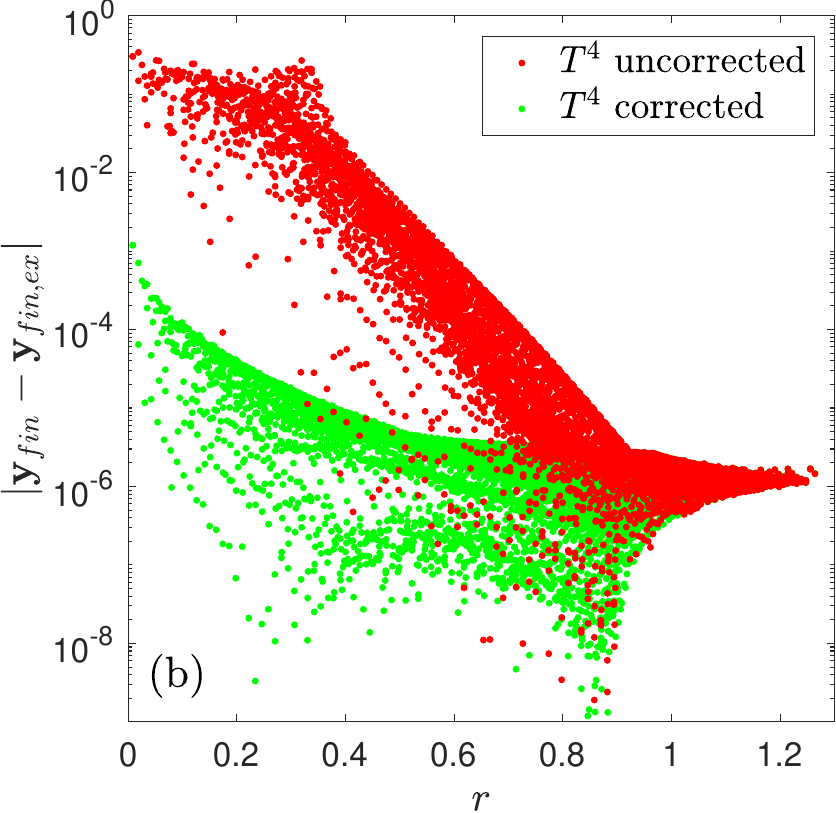}
\vskip-0.05truecm
\caption{Particles initially positioned in the square normal to $\bf{U}_{\infty}$ indicated by the green dashed line in figure 6b are evolved with the fluid velocity until they reach the parallel plane on the leeward side of the sphere indicated by the red dashed line in figure 6b. Scatter plots are shown for (a) the traversal time $T_{fin}$ and (b) the error in the final position $|\yb_{fin}-\yb_{fin,ex}|$ as a function of initial radial coordinate $r$, computed with the exact velocity, the uncorrected velocity, and the corrected velocity, as indicated.}
\label{F:1Sshadow2}
\end{figure}

The maximum error at distance $d$ is not easy to compute for fully 3D problems such as the three-ellipsoid case studied below in \S \ref{S:three-ellipsoid}.  The errors in a cross-section also do not fully capture the incurred error in a general 3D flow.  We thus present an alternative approach to measure the incurred error that will be useful below, when considering flow past several objects.  We track the motion of particles under the computed velocity field as they move from a region on the windward side of the sphere, where the flow is coming from, to the other, leeward, side, and compute the integral error incurred in the particle position and in the time it takes to traverse the region.  

Specifically, we place $75^2$ particles on a square of size $0.9\times 0.9$ normal to the incoming flow at distance 0.1 from the sphere.  This square is indicated in figure \ref{F:1Sstream}(b) by the green dashed line.  For each particle, its trajectory is computed by integrating the velocity field along particle paths using RK4 with sufficiently small timesteps to yield resolved results (here, $\Delta t=0.02$).  The trajectory is computed until the particle has reached a plane parallel to the windward square on the leeward side, at distance 0.1 from the sphere, indicated by the red dashed line in figure \ref{F:1Sstream}(b). We record the time $T_{fin}$ and position  $\yb_{fin}$ when the particles have reached the leeward plane.  The particles are evolved using the  uncorrected velocity and the corrected velocity with $n=40$, and using the exact velocity, as reference.

Figures \ref{F:1Sshadow2}(a) shows a scatter plot of the particles' traversal time $T_{fin}$ as a function of 
the initial radial coordinate $r$ (the particles initial distance to the symmetry axis) for all particles in the plane, where the trajectories are computed using the  uncorrected velocity (red dots), the corrected velocity (green dots), and the exact velocity (blue dots). Not represented here are the particles that, computed with the uncorrected velocity, never make it to the other side. Of those shown, the traversal times computed without corrections are up to 100 times larger than the actual values. The corrected results are indistinguishable from the exact values. 

Figure \ref{F:1Sshadow2}(b) shows the error in the final position  $|\yb_{fin}-\yb_{fin,ex}|$, as a function of $r$. Here, the exact values $\yb_{fin,ex}$ follow by symmetry. While the cumulative errors using the corrected velocity increase towards the center, $r=0$, where the traversal times are longer, they are less than $10^{-3}$ in all cases. The errors using the uncorrected velocity are often  100-1000 times larger.

\subsection{Flow past an ellipsoid}
We now consider a standard ellipsoid (\ref{E:ellipse}) with major axes $(a,b,c)=(3,2,1)$ in a far-field $\Ubinf=(1,0,0)$. This case adds complexity to the previous example: (i) The density is not known analytically. It and all its  derivatives are precomputed at the grid points as described in \S 2.  The values at the base point $\xb_b$ required for the coefficients $c_{pqk}$ are obtained by interpolation. (ii) The latitude-longitude grids are no longer orthogonal, so that $c_{\alf\bet}\ne0$. Thus the recursion in section \S 3.4 is necessary to compute $\int H_{pqk}$. (iii) The DLP is not zero. (iv) The exact solution is not known.  Here, the mesh used is defined by $m$ with $n_1=4m, m_1=m, n_2=3m, m_2=2m$. We use the solution with $m=320$ as a reference to compute errors shown below.
The precomputed density at all gridpoints was obtained as described in \S 2, using $p=30$, and GMRES tolerance $10^{-12}$. 

\begin{figure}[t!]
\includegraphics[trim=0 0 0 0, clip, width=0.46\textwidth]{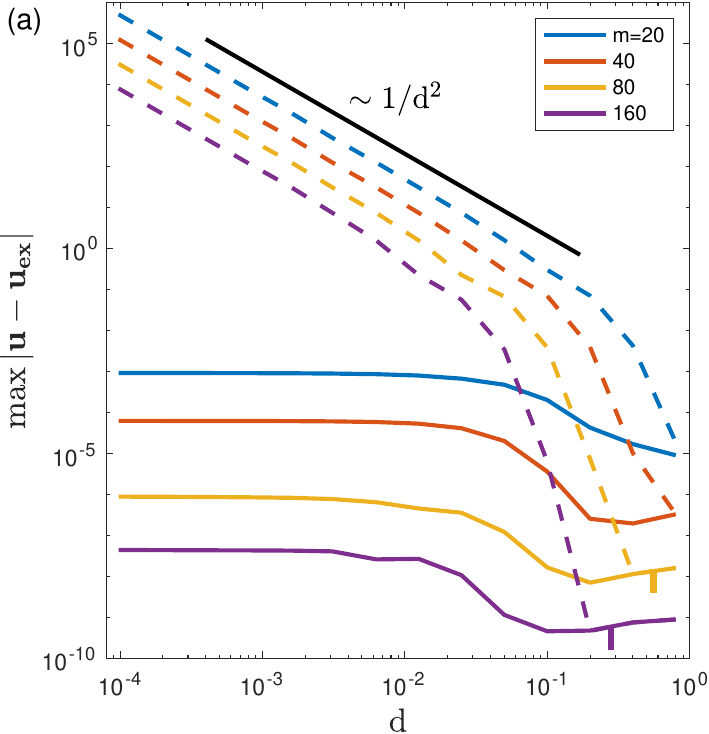}
\vskip-3.2truein\hskip3.0truein
\includegraphics[trim=30  112 95 90, clip, width=0.46\textwidth]{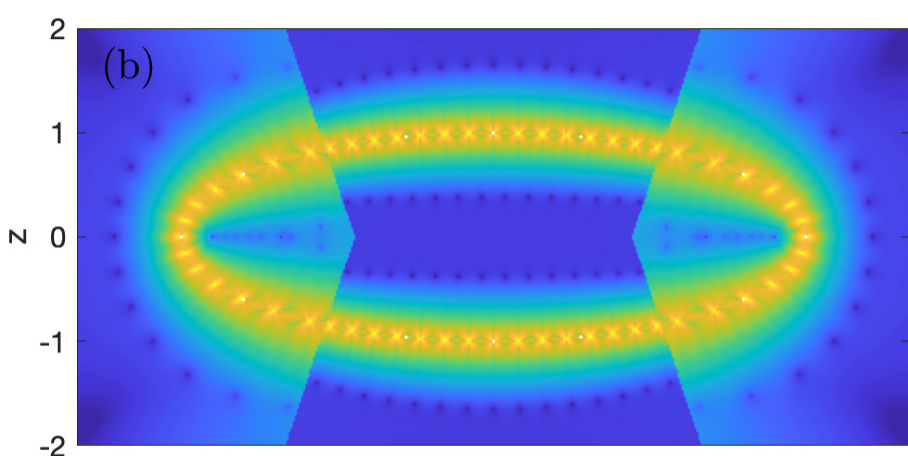}
\vskip-0.truein\hskip3.0truein
\includegraphics[trim=30  15 95 90, clip, width=0.46\textwidth]{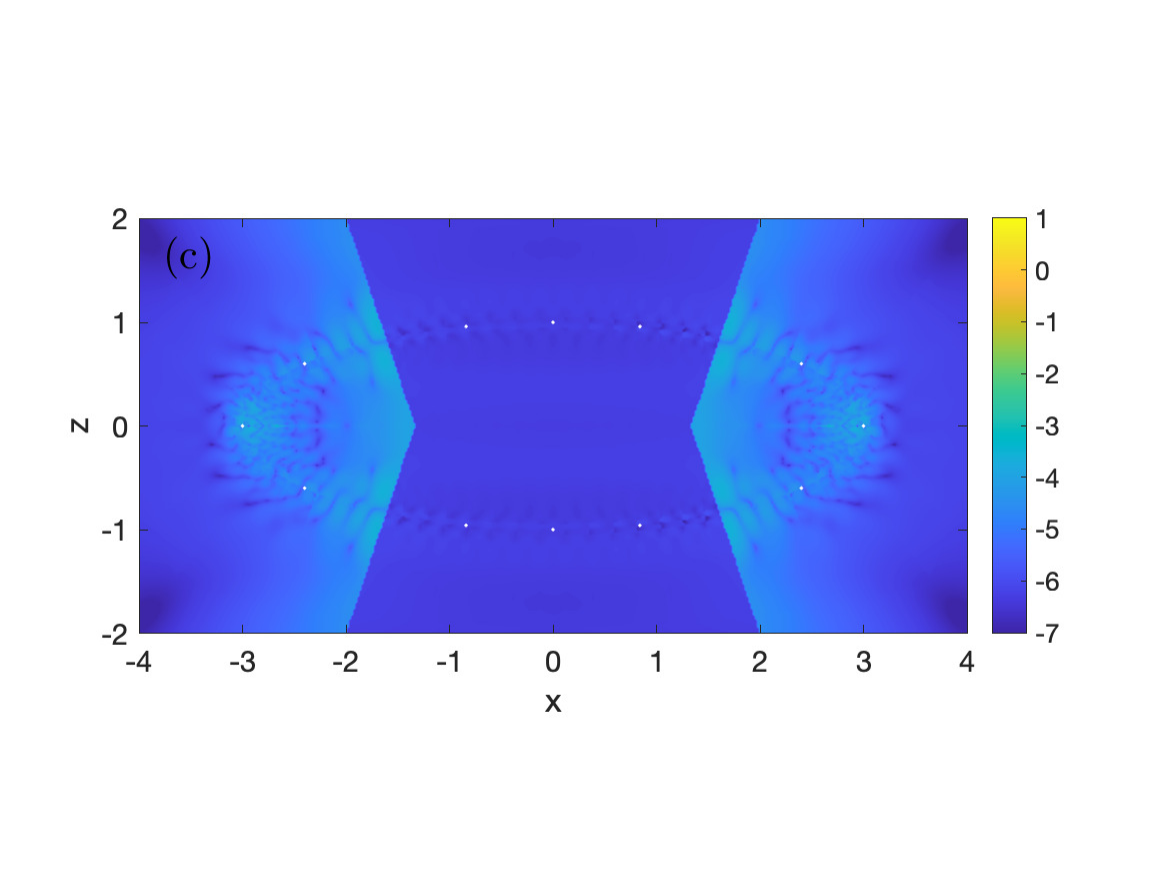}
\vskip-3.0truein\hskip6.0truein
\includegraphics[trim=470  15 40 100, clip, width=0.064\textwidth]{1E-crosserr-corr}
\vskip-2.6truein\hskip6.27truein {\scriptsize +}
\vskip2.4truein
\caption{(a) Maximal error in the velocity past the ellipse at distance $d$ from the surface, vs $d$, computed without corrections (dashed curves) and with corrections (solid curves), using the indicated values of $m$. 
For $m\ge 80$, a short vertical line indicates the distance $d_{max}=6\Delta s_{max}$ past which no corrections were applied.
(b,c) Distribution of the error in the cross section $y=0$, computed with $m=20$, without corrections in (b) and with corrections in (c).}
\label{F:1Econv}
\end{figure}

Figure \ref{F:1Econv} plots the errors in the computed flow past the ellipsoid.  Figure \ref{F:1Econv}(a) plots the maximum error in the velocity over all points at distance $d$ from the ellipsoid, computed with the indicated values of $m$, without corrections (dashed curves) and with corrections (solid curves).  The uncorrected errors satisfy the  scaling $O(h^2/d^2)$, as expected for the double layer.  With corrections, the errors are reduced by factors of $10^9$--$10^{11}$.  The corrected errors are $O(h^4)$, where $h=\pi/m$, uniformly as $d\to 0$.

Figures \ref{F:1Econv}(b,c) show the distribution of the error in the  cross-section $y=0$, computed with $m=20$, inside and outside the ellipse, without corrections in (a) and with corrections in (b), on a logarithmic scale.  Figure \ref{F:1Econv}(b) clearly shows large errors in the uncorrected velocity near the boundary.  The maximum errors are of order $10^6$, as seen in figure \ref{F:1Econv}(a), and occur near the grid points.  The target points using the two distinct grids are clearly visible, since the resolution on grid 2 in this plane is finer than on grid 1; the resolution of the two grids was chosen to match in the $z=0$ plane. Figure \ref{F:1Econv}(c) shows the  error after adding corrections. The maximum corrected errors are $10^{-3}$.

As in figure \ref{F:1Sconv}, the maximal distance $d_{max}=6\Delta s_{max}$ above which no correction is applied is indicated by a short vertical line for $m\ge 80$. As before, the figure shows that for $d\approx d_{max}$ the corrections no longer improve the standard fourth order trapezoid rule, making no further correction necessary for larger $d$.
For $m=20$ used in figures (b,c), this distance is bigger than 1, with most uncorrected points lying outside the domain shown.

\subsection{Stokes flow past 2 spheres} 

\begin{figure}
 \centering
\includegraphics[trim=90 20 100 30, clip, height=0.52\textwidth]{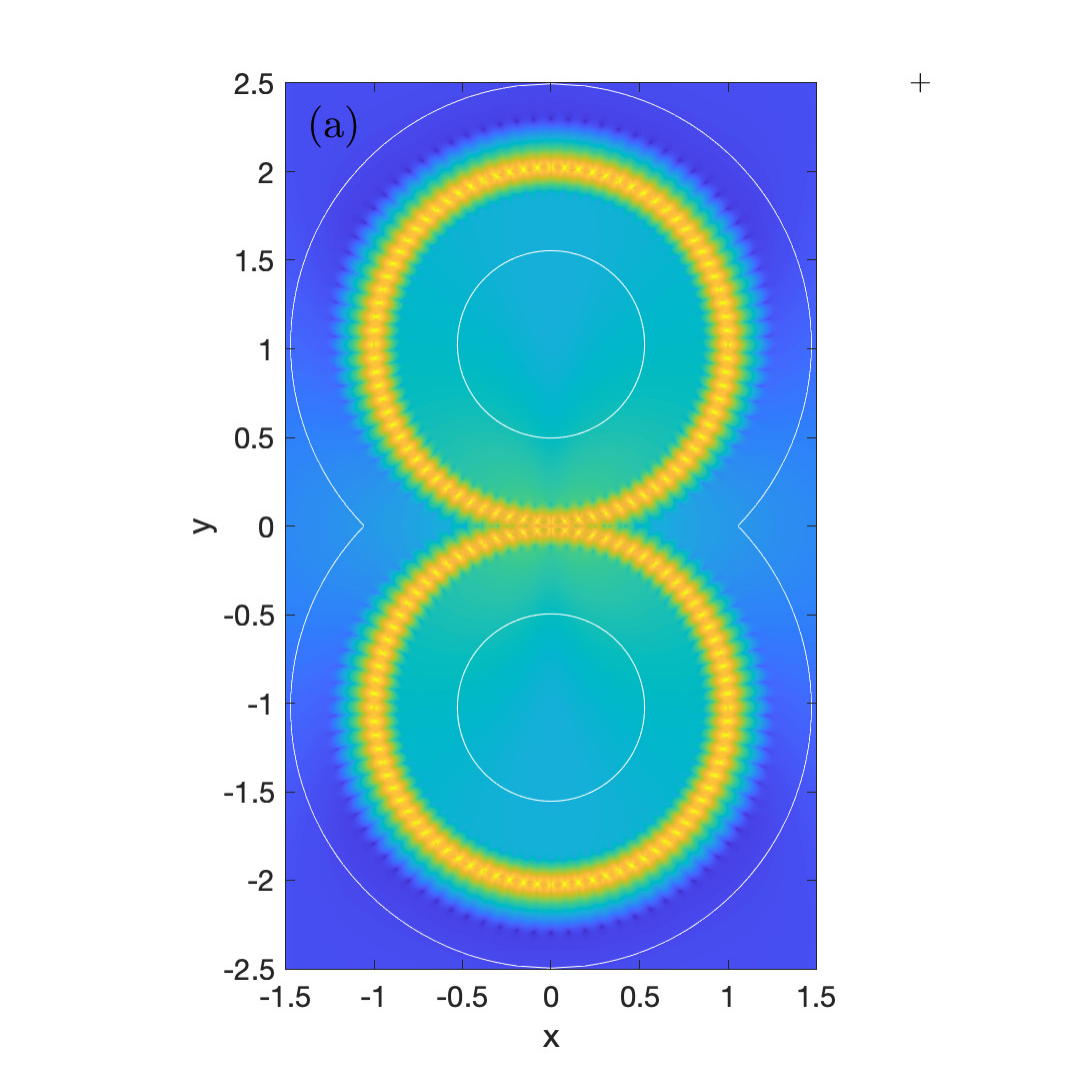}
\includegraphics[trim=115 20 100 30, clip, height=0.52\textwidth]{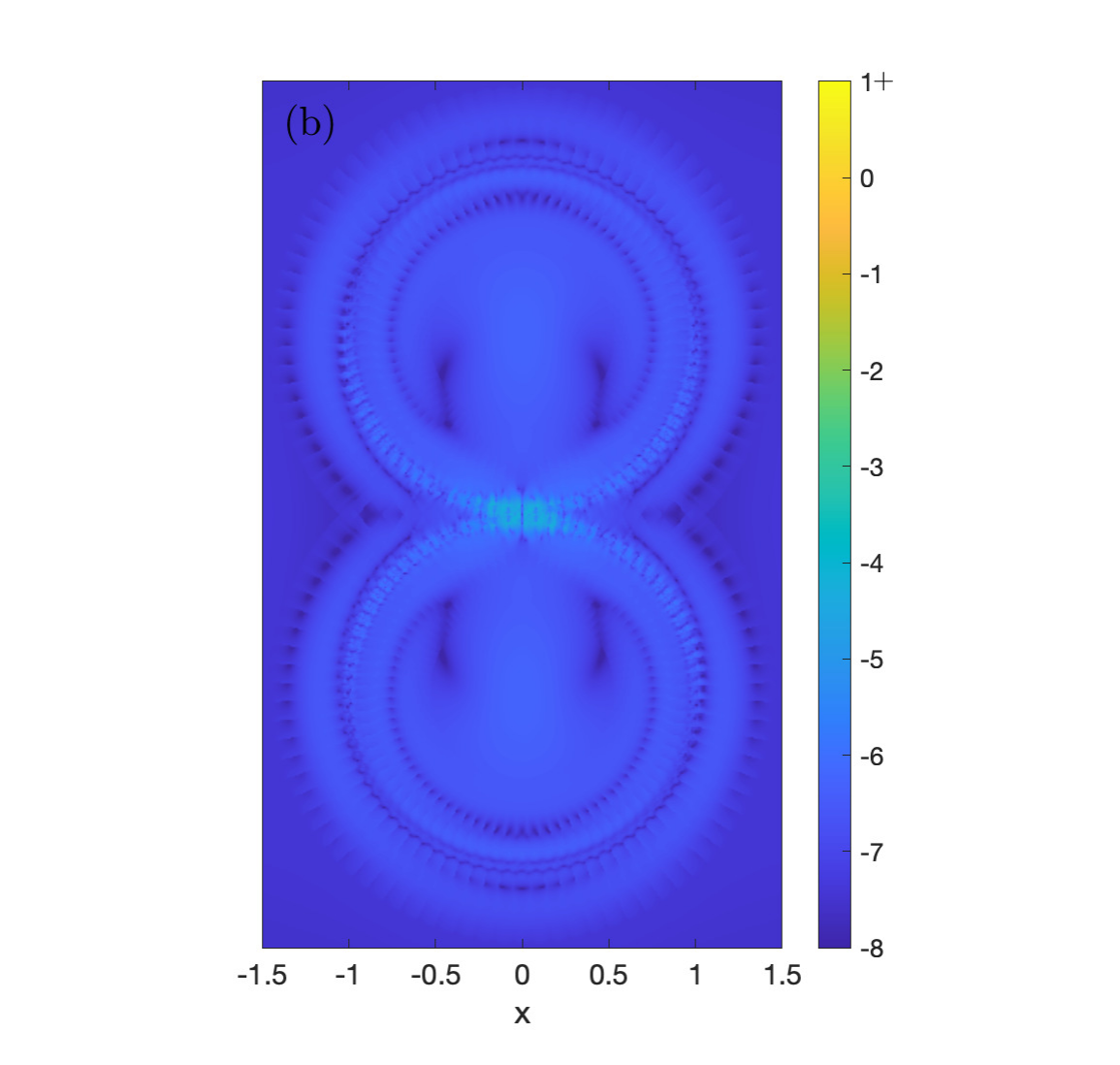}
\caption{Velocity error $|\mathbf{u}-\mathbf{u_{ex}}|$ in the plane $z=0$,
computed with $n=80$, (a) without and (b) with corrections, on a logarithmic scale. $\ub_{ex}$ is obtained with $n=320$.
{The region $d<6\Delta s$ within which corrections are applied is outlined in white.} }
\label{F:2Scross}
\end{figure}

Next, we consider flow past two spheres of radius 1, centered at $(0,\pm (1+d/2), 0)$, with far-field $\Ubinf=(1,0,0)$, with a gap of size $d=0.05$ between them. Each sphere is discretized by trapezoidal grids with
$n_1=n_2=n$ and $m_1=m_2=n/2$. Recall from \S 2 that the density on these grids is obtained by solving the linear system (2.5) to cancel the induced velocity on a spectral grid of order $p$ on each body. This system is solved iteratively with GMRES.
Here we use $p=60$ and GMRES tolerance $10^{-12}$.

Unlike in our previous examples, the velocity on the spheres contains a singular self-induced component, and a regular component induced by the second sphere, which is near-singular for target points near it.
Each iteration in GMRES requires evaluating this velocity on the spectral grid on each sphere using the current density. The singular integrals are computed on the spectral grid. The regular or near-singular integrals are computed after interpolating the density to the trapezoidal grids. The near-singular integrals must be computed with  corrections for accuracy. 

Figure 10 plots the error in the velocity obtained with $n=80$, where both the density and the velocity are computed without corrections in (a), and with corrections in (b). The error is obtained by comparison with results using $n=320$ and corrections.
The region $d\le 6\Delta s$ within which the velocity is corrected is outlined in white, in (a). Note that the errors between (a) and (b) also differ outside this region since the density in (a), obtained without corrections, is different than the density in (b), obtained with corrections. 
 As in earlier results, the corrected error in (b) is uniformly small across the boundary away from the gap between the sphere.  The error is largest in the gap, where the boundary density is larger in magnitude and less well resolved.

\subsection{Flow past several ellipsoids}
\label{S:three-ellipsoid}
The last example is that of a flow past {the three ellipsoids shown in figure \ref{F:examples}(d). They are} positioned somewhat randomly, but so that they are close to each other, with centers $(x_c,y_c,z_c)$, semi-axes lengths $a, b, c$, and  orientation angles $(\phi,\theta,\psi)$, 
 {as} given in table \ref{T:params}.
 \begin{table}[b!]
 \centering
 \renewcommand{\arraystretch}{1.3}
 \begin{tabular}{|c|c|c|c|c|c|}
 \hline
 & $a,b,c$ & $x_c,y_c,z_c$ & $\phi,\theta,\psi$ &$m_1,n_1$& $m_2,n_2$\\
 \hline
 Ellipsoid 1 & $1.80, 2.70, 0.90$ &$-0.95, ~2.0, ~0.5$ & $0, 0, 0$~& 20,80& 32,98\\
 Ellipsoid 2 & $2.50,2.50,1.25$ &$-1.0,-2.0,-0.5$ & $\frac{\pi}{3}, ~\frac{\pi}{4}, ~\frac{7\pi}{8}$~& 20,80& 40,80\\
 Ellipsoid 3 & $2.70,1.35,1.35$ &$2.0,~0.0,-0.5$ & $ \frac{\pi}{2}, -\frac{\pi}{3}, -\frac{\pi}{5}$&30,60&42,42\\
 \hline
 \end{tabular}
 \caption{Parameters describing positions and discretizations of three ellipsoids.}
 \label{T:params}
 \end{table}
The corresponding gaps between any two of these ellipsoids, defined as the smallest distance between them, range between 0.21 and 0.25, which are less than 5\% of the diameters of the ellipsoids.  Table \ref{T:params} also lists the discretization values $m_1,n_1,m_2,n_2$ used in all the results shown below. We have not attempted to optimize these values; they are chosen so that the ellipsoids have comparable mesh resolutions along their equators. The background velocity is $\Ubinf=(1,0,0)$.

\begin{figure}[!t]
 \centering
\includegraphics[trim=30 35 90 60, clip, height=0.36\textwidth]{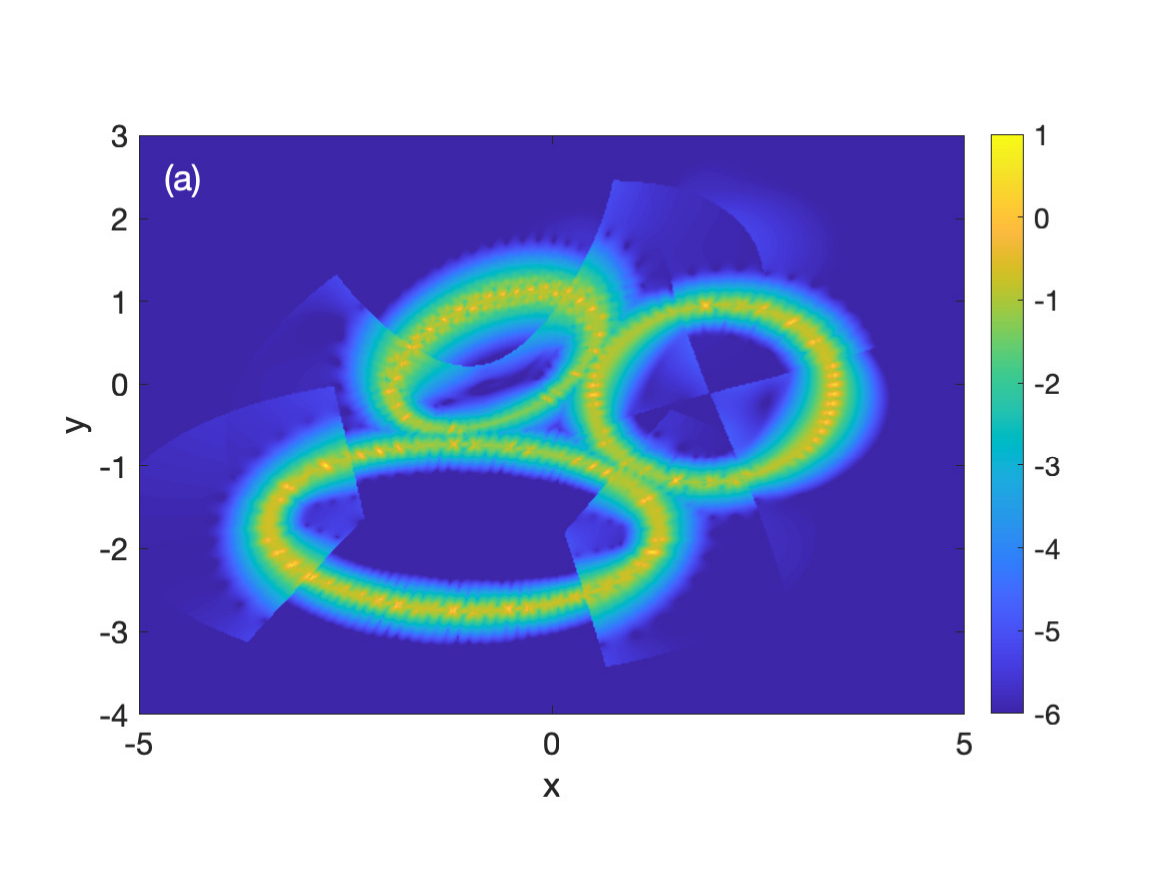}
\includegraphics[trim=62 35 50 60, clip, height=0.36\textwidth]{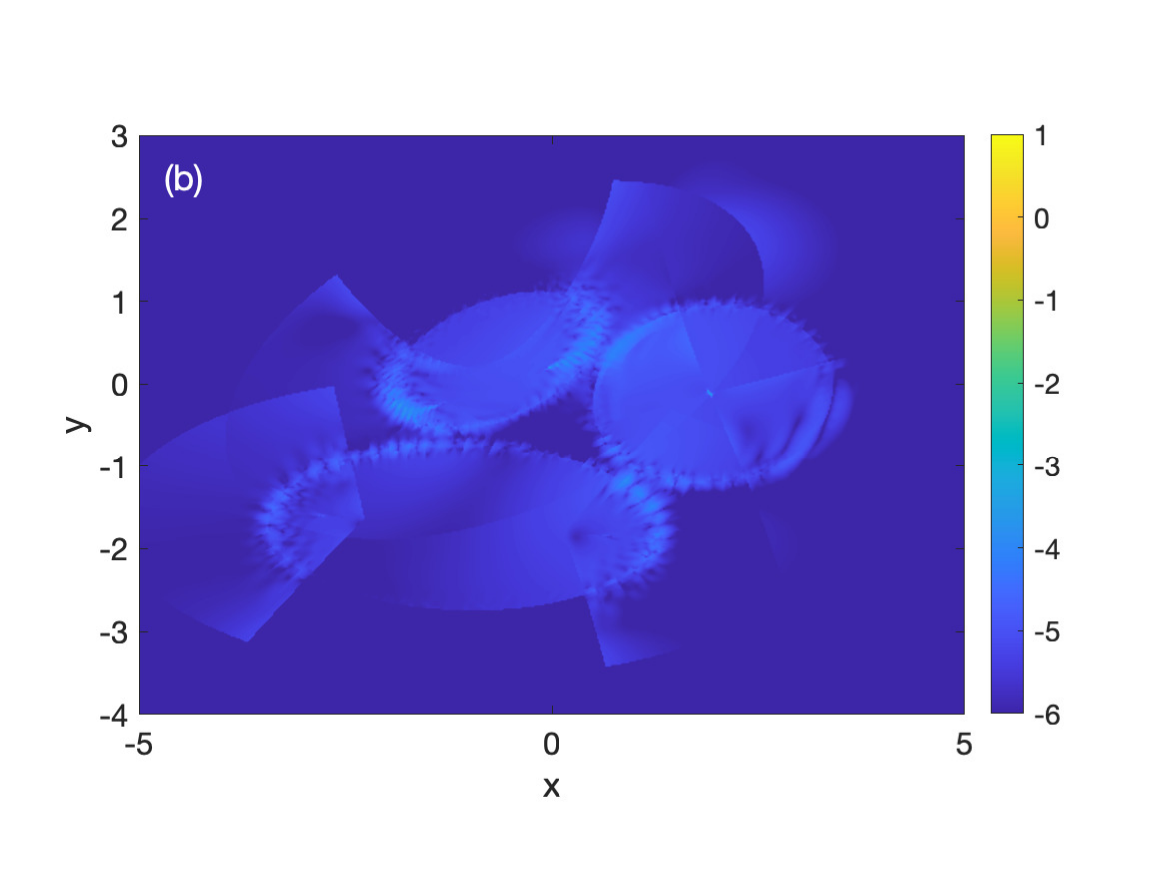}
\caption{Magnitude of error in computed velocity around 3 ellipses, in a plane through the 3 points closest to each pair of ellipses (projected onto the $x$-$y$ plane). Left: uncorrected. Right: corrected.  $m=20$.  The error is obtained by comparing with results obtained with three times the resolution.} 
\label{F:3Ecross}
\end{figure}

\begin{figure}
 \centering
\includegraphics[trim=0 0 0 0, clip, width=0.415\textwidth]{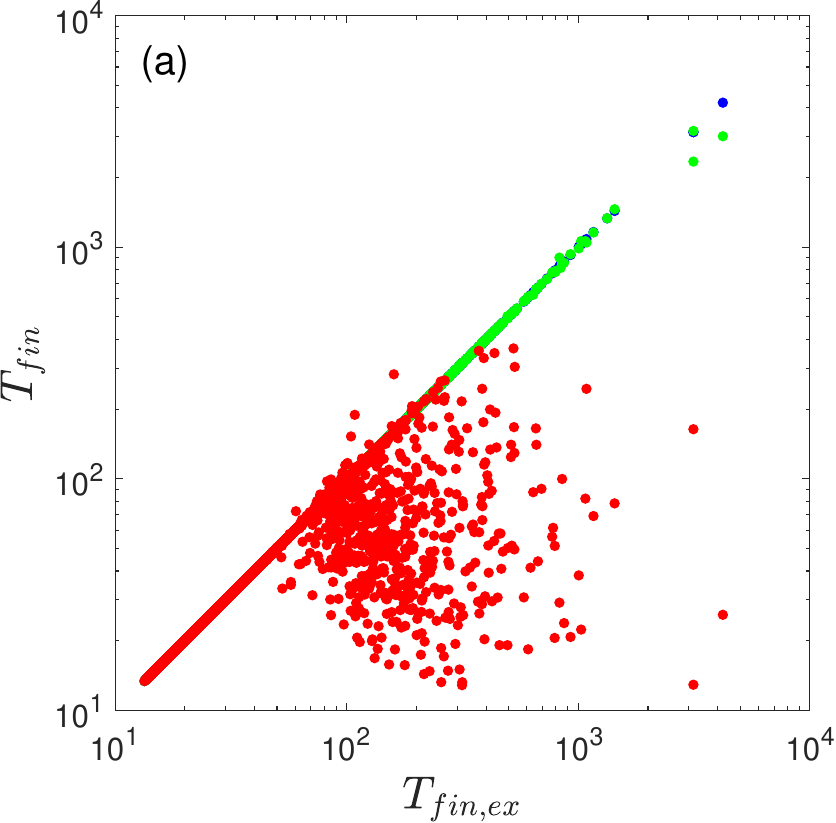}
\includegraphics[trim=0 0 0 0, clip, width=0.415\textwidth]{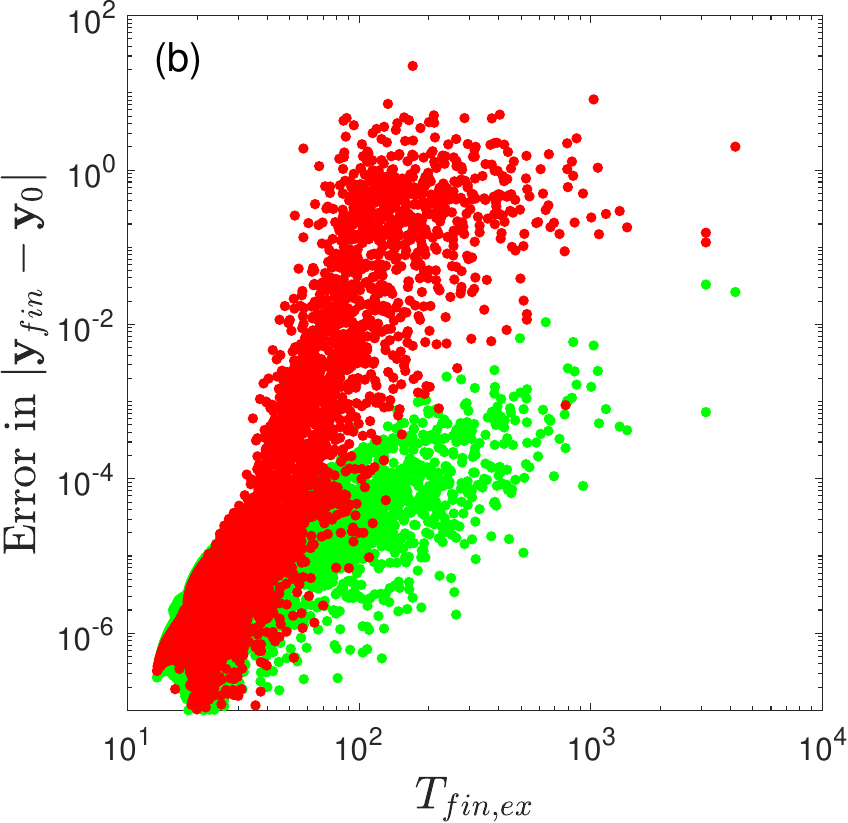}
\vskip-0.3truecm
\caption{Scatter plots of (a) $T_{fin}$ and (b) error in the displacement {$|\yb_{fin}-\yb_{fin,ex}|$}, as a function of $T_{fin,ex}$, by comparing with corrected results using a mesh three times as fine.}
\label{F:3Eshadow2}
\end{figure}

One measure of the error in the velocity field is shown in figure \ref{F:3Ecross}.
It plots the error magnitude in the plane through the three points closest to each pair of ellipsoids. In figure (a), the velocity and the densities are computed without corrections; in figure (b), they are computed with corrections. The error is obtained by comparing with results obtained with a mesh three times as fine. The uncorrected velocity in figure \ref{F:3Ecross}(a) has large errors near the ellipsoid boundary, with maximal value 72.1. The corrected results in figure (b) are significantly better resolved, with maximal error $9.7\times 10^{-4}$.

 {Figure \ref{F:3Eshadow2} presents a measure of the integral error in 3D that is similar to figure \ref{F:1Sshadow2}, where} we place a set of particles in a plane upstream of the ellipsoids, and measure errors in the time and position of those particles when they reach a parallel plane downstream of the ellipsoids. 
{We note that the time of traversal in viscous flow through a porous media is a quantity of interest in applications, as studied for example in \cite{quaife2020-transport}.} 
Here, we place $40^2$ particles in a square in the plane $x=-3.6$ with side length $6$ centered at $(y_c,z_c)=(0.35,-0.5)$. The initial particle position is given by their coordinates $\yb_0=(y_0,z_0)$ in this plane. We then compute the trajectory of those particles by integrating the velocity field 
using RK4 with a sufficiently small timestep, and record the time $T_{fin}$ and their positions $\yb_{fin}=(y_{fin},z_{fin})$ when they reach $x=3.6$.

{
Figure \ref{F:3Eshadow2}(a) shows a scatter plot of the time $T_{fin}$ for all initial position $\yb_0$,
as a function of the exact traversal time $T_{fin,ex}$, as obtained using a fine mesh. The red dots show results computed without corrections, the green dots show results computed with corrections, and the blue dots (mostly below the green dots) correspond to $T_{fin}=T_{fin,ex}$. Not represented in this figure are several particles trajectories that, when computed without corrections, never make it across due to large errors in the uncorrected velocity. The figure shows that the time to traverse the region is highly inaccurate without corrections. 
These uncorrected values are generally much smaller than the exact values, sometimes more than 100 times smaller.
With corrections, the traversal times are close to the exact values.

Figure \ref{F:3Eshadow2}(b) shows a scatter plot of the error $|\yb_{fin}-\yb_{fin,ex}|$ versus the exact traversal time $T_{fin,ex}$, for all particles, computed without corrections (red dots) and with corrections (green dots); the exact positions $\yb_{fin,ex}$ are again obtained using a fine mesh. The data show that, when computed without corrections, the error in the particle position after traversing the ellipsoids is up to $10^4$--$10^5$ times larger than the correct values.
}

\section{Conclusions}

This paper presents a quadrature method for near-singular integrals (\ref{E:genform}) that arise in the Boundary Integral Equation method in 3D, for target points $\xo$ close to the surface of integration. As a direct generalization of the quadrature of \cite{nitsche2021tcfd} from line integrals to surface integrals, we derive a high-order quadrature by error correction of the double Trapezoidal rule. 
The near-singular integrands are approximated by an expansion using basis functions that can be integrated to high accuracy with recursion formulas. Fourth-order corrected quadratures are derived for the Stokes kernels and their convergence is validated and analyzed at target points near a sphere or an ellipsoid. For any fixed discretization, uniform accuracy is achieved for all close targets regardless of their distance $d$ to the surface. Proper care is taken to resolve the loss of significance due to roundoff effects as $d\to0$. The quadratures are applied to simulate viscous flows past multiple objects (two spheres or three ellipsoids), where each pair of objects is separated by a small distance that is 2.5\% - 5\% of the object diameters.

To develop more efficient integral equation solvers, our correction scheme could also be used together with standard fast summation algorithms, such as the Fast Multipole Method (FMM) \cite{rokhlin1987,tornberg2008fast}; although our correction scheme uses two different grids, it is sufficient to use only one grid for all the far-field targets outside the region of correction, thus all of the far interactions can be accelerated by the FMM. In addition, a new singular quadrature scheme that uses the same two trapezoidal grids on the sphere is also under development, where the singular correction will be based on the ``zeta quadrature'' of Wu and Martinsson \cite{wu2021corrected}. In the future, this approach will allow us to use the same double-trapezoidal grids for both singular and near-singular quadrature, further simplifying the quadrature scheme.

Our quadrature correction scheme can also be extended to accommodate a greater variety of problems, which may involve different integral operators and parametric surfaces of more general shapes.
Corrections for the Laplace and Helmholtz kernels, with applications to vortex flows and wave scattering problems, are currently being explored. 
For surfaces other than ellipsoids, our correction scheme can be applied almost without modification to deformable particles in Stokes flow \cite{shravan-vesicle3d}, and can be extended to even more general surfaces by combining with domain decomposition methods based on the partition-of-unity approach \cite{bruno2001fast,beale2016} where a surface is represented as a collection of overlapping quadrilateral patches; these developments will be reported in the future.

\appendix
\section{Recursion for $ F_{pk}, G_{qk}$}
\label{app_A}
The integrals
\begin{equation}
F_{pk}=\int_a^b \frac{u^p}{\rho_u^{2k+1}}\,du~,\quad
G_{qk}=\int_c^d \frac{v^q}{\rho_u^{2k+1}}\,dv~,\quad
\end{equation}
are evaluated recursively, as follows. 
Set
\begin{subequations}
\begin{align}
&F_{00}=-\ln(-u-{C}v+\rho_u)\Big|_{u=a}^{u=b}\\ 
&F_{0k}= \frac{1}{[(1-{C}^2)v^2+1](2k-1)}\Big[\frac{u+{C}v}{\rho_u^{2k-1}}+2(k-1)F_{0(k-1)}\Big],~k=1:\kmax
\\
&F_{p0}= \frac{1}{p}\Big[u^{p-1}\rho_u - (2p-1){C}vF_{(p-1)0}-(p-1)(v^2+1)F_{(p-2)0}\Big] ,~p=1:\mmax\\
&F_{pk}= \frac{1}{2k-1} \Big[-\frac{u^{p-1}}{\rho_u^{2k-1}} +(p-1) F_{(p-2)(k-1)} \Big] -{C}v F_{(p-1)k},~ k=1:\kmax, ~p=1:\mmax
\end{align}
\end{subequations}
where all $F_{pk}$ with negative indeces are set to zero.
The values of $G_{qk}$ are computed similarly, by simply replacing $F, u, p, a, b$ respectively by $G, v, q, c, d$ above.

\section{Positive definiteness of $\calf^2\alpha^2+2c_{\alf\bet}\alpha\beta+\cbet^2\beta^2$}
\label{app_B}
For given target point $\xo$ and its orthogonal projection $\xb_b$ onto an ellipsoid $\xb(\alf,\bet)$, define the quadratic form
\begin{equation}
Q(\alpha,\beta):=A\alpha^2+2C\alpha\beta+B\beta^2
\label{E:quadratic_forms}
\end{equation}
 where the coefficients $A,B,C$ 
 are given by (using the notation of \eqref{E:calf})
\begin{eqnarray}
A&=&\calf^2= (\xb_b-\xo)\cdot\xb_{\alf\alf}+\xalf\cdot\xalf,\nonumber\\
C&=&c_{\alf\bet}= (\xb_b-\xo)\cdot\xb_{\alf\bet}+\xalf\cdot\xbet,\\
B&=&\cbet^2= (\xb_b-\xo)\cdot\xb_{\bet\bet}+\xbet\cdot\xbet,\nonumber
\label{E:calfagain}\end{eqnarray}
The quadratic form is \textit{positive definite} if $Q(\alf,\bet)> 0$ for all $(\alf,\bet)\ne(0,0)$, or, equivalently, if $A>0$, $B>0$, and $C^2<AB$.

Let $\nb$ be the outward unit normal at $\xb_b$, and $d$ be the signed distance from $\xo$ to $\xb_b$,
\begin{equation}
\xo-\xb_b = d\nb,
\end{equation}
In this section we show that the quadratic form is positive definite provided $\xo$ is either outside the ellipse, or if it is inside, at a distance $|d|$ smaller than the radius of the osculating sphere.
The osculating sphere at $\xb$ is the largest sphere inside the ellipsoid that only intersect the ellipsoid at a single point $\xb$.
Under these conditions, $A>0$ and $B>0$, as shown in \cite{nitsche2021tcfd}. It remains to show that $C^2<AB$.

We first define the following.
\begin{enumerate}
\item The coefficients of the first and second fundamental forms at $\xb_b$ are
\begin{align*}
E&:=\xb_\alpha\cdot\xb_\alpha,~~F:=\xb_\alpha\cdot\xb_\beta,~~G:=\xb_\beta\cdot\xb_\beta\\
L&:=\nb\cdot\xb_{\alpha\alpha},~~M:=\nb\cdot\xb_{\alpha\beta},~~N:=\nb\cdot\xb_{\beta\beta}\,.
\end{align*}
\item The mean curvature $H$ and the Gaussian curvature $K$ at $\xb_b$ are given by
$$H := \frac{EN-2FM+GL}{2(EG-F^2)}\,~~ K:=\frac{LN-M^2}{EG-F^2},$$
\item The principal curvatures at $\xb_b$ are 
$$\kappa_1:= H-\sqrt{H^2-K}~~\text{and}~~\kappa_2:=H+\sqrt{H^2-K}.$$
\end{enumerate}
The coefficients \eqref{E:calfagain} can be written as
\begin{eqnarray}
A&=& E-Ld,\nonumber\\
C&=& F-Md,\\
B&=&  G-Nd.\nonumber
\label{E:calf_via_fundamental_form}
\end{eqnarray}
Then the positive definiteness condition $AB>C^2$ can be equivalently written as
$$1-2dH+d^2K>0.$$
Solving the inequality above gives
\begin{equation}
\frac{1}{d}<\kappa_1~~~\text{or}~~~\frac{1}{d}>\kappa_2.
\label{E:d_kappa}
\end{equation}
Note that $d>0$ when $\xo$ is outside of the ellipsoid, and $d<0$ when $\xo$ is inside. Also, $\kappa_1<\kappa_2<0$
 at any $\xb_b$ on an ellipsoid. 
 
 For any point $\xo$ outside the ellipsoid we have $d>0$, thus \eqref{E:d_kappa} is satisfied and the quadratic form \eqref{E:quadratic_forms} is positive definite. 
On the other hand, if  $\xo$  is inside the ellipsoid, then $d<0$, in which case \eqref{E:d_kappa} is satisfied when $|d|<1/|\kappa_1|$. In that case also $A,B>0$. 
Equation \eqref{E:d_kappa} is also satisfied when $|d|>1/|\kappa_2|$, but in that case, $A$ and $B$ are not positive.

Note that $R_{osc}:=1/|\kappa_1|$ is the radius of the osculating sphere along the tangential direction of largest negative curvature at $\xb_b$. It follows that the quadratic form \eqref{E:quadratic_forms} is positive definite when either $\xo$ is outside the ellipsoid, or the distance of $\xo$ from inside is smaller than $R_{osc}$.
\bibliographystyle{apalike}
\bibliography{}
\end{document}